\documentclass[letterpaper, 10 pt, conference]{ieeeconf}

\bstctlcite{bstctl:etal, bstctl:nodash, bstctl:simpurl}

\usepackage{anyfontsize}
\usepackage{array}
\usepackage[normalem]{ulem}
\usepackage{stfloats}
\usepackage{mathtools}
\usepackage[utf8]{inputenc}
\usepackage{tabularx}
\usepackage{booktabs}
\usepackage{mathrsfs}
\usepackage{graphics}

\usepackage{amsmath}
\usepackage{amssymb}
\usepackage{amsthm}
\usepackage{cite}
\usepackage{bm}
\usepackage{acronym}
\usepackage{paralist}
\usepackage{float}
\usepackage[dvipsnames]{xcolor}
\usepackage{epstopdf}
\usepackage{multicol}
\usepackage{tikz}
\usepackage{hyperref}
\usepackage{graphicx}
\usepackage{soul}
\usepackage{macros}

\fnbelowfloat

\makeatletter
\hypersetup{colorlinks=true}
\AtBeginDocument{\@ifpackageloaded{hyperref}
  {\def\@linkcolor{blue}
  \def\@anchorcolor{red}
  \def\@citecolor{red}
  \def\@filecolor{red}
  \def\@urlcolor{black}
  \def\@menucolor{red}
  \def\@pagecolor{red}
\begingroup
  \@makeother\`%
  \@makeother\=%
  \edef\x{%
    \edef\noexpand\x{%
      \endgroup
      \noexpand\toks@{%
        \catcode 96=\noexpand\the\catcode`\noexpand\`\relax
        \catcode 61=\noexpand\the\catcode`\noexpand\=\relax
      }%
    }%
    \noexpand\x
  }%
\x
\@makeother\`
\@makeother\=
}{} }
\makeatother

\newtheorem{AssumptionEnv}{Assumption}

\newtheorem{ProblemEnv}{Problem}
\newenvironment{Problem}[1][]
{\begin{ProblemEnv}}
{\hfill$\diamond$\end{ProblemEnv}}

\newtheorem{RemarkEnv}{Remark}
\newenvironment{Remark}[1][]
{\begin{RemarkEnv}}
{\hfill$\blacklozenge$\end{RemarkEnv}}

\newtheorem{TheoremEnv}{Theorem}

\newtheorem{LemmaEnv}{Lemma}

\newtheorem{PropositionEnv}{Proposition}
\newenvironment{Proposition}[1][]
{\begin{PropositionEnv}}
{\hfill$\vartriangle$\end{PropositionEnv}}

\newtheorem{CorollaryEnv}{Corollary}

\newtheorem{DefinitionEnv}{Definition}
\newenvironment{Definition}[1][]
{\begin{DefinitionEnv}}
{\hfill$\bullet$\end{DefinitionEnv}}

\DeclareMathOperator{\trace}{\mathrm{tr}}

\newcommand{\bequ}{\begin{eqnarray}}
\newcommand{\eequ}{\end{eqnarray}}

\newcommand{\HH}{\mathcal{H}}
\newcommand{\VV}{\mathcal{V}}

\newcommand{\Ito}{It\^{o}}

\IEEEoverridecommandlockouts

\def\BibTeX{{\rm B\kern-.05em{\sc i\kern-.025em b}\kern-.08em
    T\kern-.1667em\lower.7ex\hbox{E}\kern-.125emX}}

\definecolor{AA}{RGB}{34,139,34}

\newcommand{\Sph}{\mathbb S^2}
\newcommand{\Grnd}{\mathbb{G}_\mathrm{rnd}}
\newcommand{\Q}{\mathsf Q}
\newcommand{\TQ}{\mathsf{T}\mathsf Q}

\newcommand{\piQ}{\pi}
\newcommand{\LC}{\nabla}
\newcommand{\Sasaki}{\mathbb{G}^{\mathrm S}}
\newcommand{\LCs}{\overline{\nabla}}
\newcommand{\Riem}{R}
\newcommand{\Hspace}{\mathsf H}
\newcommand{\Vspace}{\mathsf V}
\newcommand{\Hor}{\mathrm H}
\newcommand{\Ver}{\mathrm V}
\newcommand{\dd}{\mathrm d}
\newcommand{\cD}{\mathrm D}
\newcommand{\Hess}{\operatorname{Hess}}
\newcommand{\muS}{\mu^{\mathrm S}}
\newcommand{\divS}{\operatorname{div}_{\mu^{\mathrm S}}}
\newcommand{\gradS}{\operatorname{grad}^{\mathrm S}}

\IEEEoverridecommandlockouts

\begin{document}

\bstctlcite{IEEE_b:BSTcontrol}

\title{\LARGE \bf
Geometric Stochastic Dynamics and Fokker–Planck Equation on \(\mathsf{T}\mathbb S^2\)
\thanks{This research has been supported by NJIT's startup funds.}
\thanks{H. Mahmood and A. Akhtar are with the
Dept. of Mechanical \& Industrial Engineering at the New Jersey Institute of Technology, Newark, NJ 07102, USA (email: \texttt{\{hm576, adeel.akhtar\}@njit.edu}).}}

\author{Hamza Mahmood \and Adeel Akhtar}
\maketitle

\begin{abstract}
This paper develops a geometric (coordinate-free) stochastic model of a controlled point mass evolving on the unit two-sphere $\mathbb S^2$, together with the Fokker-Planck equation governing the probability density of the resulting stochastic process. The state manifold is the tangent bundle $\mathsf{T}\mathbb S^2$. The Levi-Civita connection of the round metric is used to construct the horizontal-vertical splitting of $\mathsf{T}(\mathsf{T}\mathbb S^2)$ and to convert the intrinsic second-order dynamics 
of the system into a first-order controlled vector field on $\mathsf{T}\mathbb S^2$. Endowing the tangent bundle with the Sasaki metric then permits stochastic perturbations to be introduced intrinsically and independently in horizontal, i.e., position, and vertical, i.e., velocity, directions. 
The corresponding infinitesimal generator is obtained in horizontal-vertical Hessian form. 
The probability density is defined with respect to the Sasaki Riemannian volume measure. The Fokker-Planck equation is derived intrinsically in both Stratonovich and It\^o forms. In the It\^o representation, a globally defined horizontal-vertical diffusion tensor yields a compact double-divergence Fokker-Planck PDE.
No coordinates or embedded Euclidean representation of the sphere dynamics are used in the development. Numerical simulations illustrate the stochastic motion of a controlled particle ensemble on $\mathbb S^2$ toward an antipodal target configuration. 
\end{abstract}

\section{Introduction}
\label{sec:Introduction}
Stochastic perturbations arise naturally in mechanical systems evolving on nonlinear configuration spaces, where uncertainty in motion, actuation, and environmental interactions must be described without violating the underlying geometric constraints~\cite{Mahmood2026, lee2013stochastic, chirikjian2011}. When the configuration space is a Riemannian manifold, the formulation of stochastic differential equations (SDEs) requires particular care, since stochastic increments and their associated drift corrections must be compatible with the manifold geometry~\cite{LeeChirikjian2025}. Beyond individual stochastic trajectories, the Fokker--Planck equation provides a deterministic description of how the probability density of the state evolves over time~\cite{yosida1949integration, ito1950stochastic, ito1953stochastic}. A geometrically consistent formulation of both the SDE and its forward equation is therefore important for uncertainty propagation, probabilistic analysis, and stochastic control on manifolds \cite{Hsu2002,LeeChirikjian2025}.

\begin{figure}
    \centering
    \includegraphics[scale=0.5]{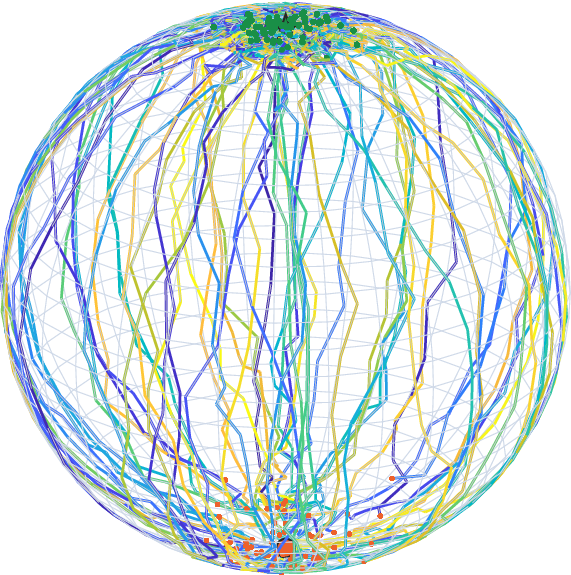}
\caption{Sample paths of an ensemble of $N = 70$ controlled
stochastic point masses on $\mathbb S^2$. Each curve
represents the configuration-space projection of a
stochastic trajectory on $\mathsf{T}\mathbb S^2$. The initial and target configurations
are the antipodal points $q_{\mathrm I}=(0,0,-1)$ and
$q_\star=(0,0,1)$, respectively.}
\label{fig:sample_path_trajectories}
\end{figure}

The geometric formulation of mechanical control systems has
been systematically developed using affine connections
and the intrinsic description of second-order dynamics
\cite{BulloLewis2005}, primarily within a deterministic
modeling and control framework. By comparison, stochastic extensions of such geometric control formulations have received comparatively less attention~\cite{mahmood2026geometricbridge, mahmood2026equivalencestochasticcontrolpath, LeeChirikjian2025}. In this work, we consider a fully actuated controlled
point mass on $\mathsf{T}\Sph$ and incorporate
independently adjustable horizontal and vertical
stochastic perturbations. We also derive the associated
intrinsic Fokker--Planck equation with respect to
the Sasaki Riemannian volume measure. 

In the literature,
stochastic systems are formulated through tangent-bundle
constructions and manifold-valued Brownian motion in
\cite{Duncan1979}, although the treatment does not
specialize to independently prescribed horizontal and
vertical diffusion on the Sasaki tangent bundle. Geometric
Langevin equations on embedded submanifolds, including
spherical configurations and constant-speed variants,
are developed in \cite{GrothausStilgenbauer2013};
these models primarily introduce stochastic forcing
through the velocity dynamics rather than independent
horizontal disturbances of the configuration.
Damped, stochastically forced geodesic motion on
Riemannian manifolds is investigated in
\cite{BirrellEtAl2017}, with emphasis on the small-mass
limit and the resulting noise-induced drift, rather
than the forward density evolution of a controlled
finite-inertia system. Kinetic Brownian motion is
studied in \cite{AngstBailleulTardif2015} on the unit
tangent bundle, where the velocity magnitude remains
fixed, in contrast to dynamics on the full tangent
bundle. Related studies examine random perturbations
of geodesic equations and their Brownian-motion limits
\cite{Li2016}, as well as stochastic averaging of
affine-connection control systems driven by rapid
random vibrations \cite{SuttnerEbenbauer2026};
their principal objectives differ from deriving a
forward equation for a controlled diffusion with
independent position and velocity noise. More generally,
Langevin-type dynamics on abstract smooth manifolds
are analyzed in \cite{GrothausMertin2022}, with
particular emphasis on hypocoercivity and exponential
convergence to equilibrium. Recent work further
develops kinetic Fokker--Planck theory on tangent
bundles in curved geometry \cite{Villani2026},
focusing on velocity-space diffusion, hypoelliptic
regularization, and equilibration. In this paper, we present geometric and stochastic frameworks to develop a unified intrinsic description of the controlled stochastic dynamics and probability density evolution of a fully actuated point mass on $\Sph$, with state manifold $\mathsf{T}\Sph$.

More recently, a global tensor-field formulation of the Fokker-Planck equation on Riemannian manifolds was developed in \cite{LeeChirikjianFP2026}. However, for such a formulation on a tangent bundle, we would like to make use of its richer geometric structure. To this end, the differential geometry of tangent bundles, including the horizontal-vertical decomposition and the Sasaki metric, provides a natural framework for treating configuration and velocity as components of a single manifold-valued state \cite{Sasaki1958,YanoIshihara1973}.

In this paper, we exploit the differential geometry of tangent bundles for a fully actuated point mass whose configuration evolves on the round two-sphere $\Sph$ and whose state belongs to the four-dimensional tangent bundle $\mathsf{T}\Sph$. Working on this state manifold is not a direct extension of a configuration-level SDE: a tangent vector to $\mathsf{T}\Sph$ cannot be identified canonically with a pair of vectors in $\mathsf{T}_q\Sph$, and stochastic position disturbances must be distinguished from stochastic changes of the velocity state. The Levi--Civita connection supplies the required horizontal and vertical lifts, while the Sasaki metric equips $\mathsf{T}\Sph$ with a natural Riemannian structure and volume measure. We use these constructions to obtain intrinsic Stratonovich and It\^o models with independently adjustable horizontal and vertical noise intensities and to specialize the Riemannian Fokker--Planck framework to the resulting controlled diffusion, including its horizontal--vertical diffusion tensor and equal-intensity Sasaki Laplace--Beltrami limit. The resulting description supports probability propagation and uncertainty-aware analysis of controlled motion on spherical configuration spaces and provides a basis for subsequent distribution-level control and stochastic planning applications.

The main contributions of this paper are summarized as follows:
\begin{enumerate}
    \item \textbf{Intrinsic stochastic dynamics:}
    We construct the controlled Stratonovich SDE
    \eqref{eq:stochastic-state-simplified} and its intrinsic
    It\^o representation \eqref{eq:intrinsic-Ito-state}
    on $\mathsf{T}\Sph$, with independent isotropic
    horizontal and vertical noise.

    \item \textbf{Geometric infinitesimal generator:}
    We derive the horizontal--vertical Sasaki-Hessian
    representation \eqref{eq:generator-HV-Hessian}
    and establish its equal-intensity case
    \eqref{eq:generator-equal-noise} to the
    Sasaki Laplace--Beltrami operator.

    \item \textbf{Intrinsic Fokker--Planck equation:}
    We construct the global diffusion tensor
    \eqref{eq:diffusion-tensor-HV} and derive the
    double-divergence forward equation
    (Theorem~\ref{thm:main-FP}), including its
    equal-intensity formulation \eqref{eq:FP-isotropic}.
\end{enumerate}

We emphasize that the global tensor-field formulation of the Fokker--Planck equation developed in \cite{LeeChirikjianFP2026} applies to general Riemannian manifolds and is therefore broader in scope than the present setting. However, by exploiting the Levi--Civita horizontal--vertical splitting and the Sasaki geometry of $\mathsf{T}\Sph$, we obtain a particularly simple and explicit representation of the diffusion tensor and the associated Fokker--Planck equation. 

\section{Preliminaries}
\label{sec:Preliminaries}
In this section, we recall the fundamental concepts and constructions from differential geometry required for the intrinsic formulation of controlled dynamics of a point mass on $\mathbb{S}^2$. We first introduce the configuration and state manifolds, together with the associated Levi-Civita connection. We then present the horizontal-vertical decomposition of the tangent bundle $\mathsf{T}\mathbb{S}^2$, including the connection map and the horizontal and vertical lifts. Finally, we use these constructions to express the (deterministic) fully actuated second-order system dynamics as a first-order control system on $\mathsf{T}\mathbb{S}^2$.

\subsection{Geometric Setting}
\label{subsec:geometric_setting}
Let
$
    \mathsf{Q}=\Sph
$
be the unit two-sphere equipped with its standard round Riemannian metric $\Grnd$. 
Throughout the paper, the Riemannian manifold $\Q$ is treated intrinsically; in particular, the deterministic and stochastic dynamics are formulated without embedding $\Q$ in an ambient Euclidean space.

For each $q\in\Q$, let $\mathsf{T}_q\Q$ denote the tangent space at $q$. The tangent bundle is
\begin{equation*}
    \TQ=\bigsqcup_{q\in\Q}\mathsf{T}_q\Q,
\end{equation*}
and its canonical bundle projection is
\begin{equation*}
    \piQ:\TQ\longrightarrow\Q,
    \qquad
    \piQ(v_q)=q.
\end{equation*}
The round metric $\Grnd$ determines a unique Levi--Civita connection, denoted by $\LC$. If $q(t)$ is a smooth curve on $\Q$ and $Y(t)\in \mathsf{T}_{q(t)}\Q$ is a smooth vector field along $q(t)$, its covariant derivative along $q(t)$ is denoted by
\begin{equation*}
    \frac{\cD Y}{\dd t}
    :=\LC_{\dot q(t)}Y.
\end{equation*}
The Riemann curvature tensor of $\LC$ is denoted by $\Riem$ and is taken with the convention
\begin{equation*}
    \Riem(X,Y)Z
    =\LC_X\LC_YZ-\LC_Y\LC_XZ-\LC_{[X,Y]}Z.
\end{equation*}
In particular, $\Riem(X,X)=0$.

The state of the point-mass system on $\Q$ is denoted by 
$
    x=v_q\in\TQ,
$
where $v_q$ is the velocity of the point mass at the configuration $q\in\Q$. Hence the state manifold is the four-dimensional manifold
$
    \TQ=\mathsf{T}\Sph.
$
A central geometric point is that a vector in $\mathsf{T}_{v_q}\TQ$ is not canonically an ordered pair of tangent vectors in $\mathsf{T}_q\Q$. The Levi-Civita connection supplies the additional structure required to split $\mathsf{T}_{v_q}\TQ$ into `horizontal' and `vertical' parts.

\subsection{Horizontal-Vertical Splitting of
\texorpdfstring{$\mathsf{T}_{v_q}\mathsf{T}\mathbb{S}^2$}{Tvq TS²}}

\subsubsection{Vertical subspace and vertical lift}
Fix $v_q\in\TQ$. The differential of the bundle projection is
\begin{equation*}
    T_{v_q}\piQ:\mathsf{T}_{v_q}\TQ\longrightarrow \mathsf{T}_q\Q.
\end{equation*}

\begin{Definition}[Vertical subspace]
The vertical subspace at $v_q\in\TQ$ is defined as
\begin{equation*}
    \Vspace_{v_q}\TQ
    :=\ker\!\left(T_{v_q}\piQ\right).
\end{equation*}
Thus a vector in $\Vspace_{v_q}\TQ$ changes the fiber variable $v_q$ inside the fixed vector space $\mathsf{T}_q\Q$ while leaving the base point $q$ unchanged to first order.
\end{Definition}

Because each fiber $\mathsf{T}_q\Q$ is itself a vector space, it is canonically isomorphic to the corresponding vertical space.

\begin{Definition}[Vertical lift]
For $w\in \mathsf{T}_q\Q$, the vertical lift of $w$ at $v_q$ is
\begin{equation*}
    w^{\Ver}_{v_q}
    :=\left.\frac{\dd}{\dd s}\right|_{s=0}(v_q+s w)
    \in \mathsf{T}_{v_q}\TQ.
\end{equation*}
\end{Definition}
Since the base point of the curve $s\mapsto(v_q+s w)$ is fixed,
\begin{equation*}
    T\piQ\!\left(w^{\Ver}_{v_q}\right)=0,
\end{equation*}
so $w^{\Ver}_{v_q}\in\Vspace_{v_q}\TQ$. The map
\begin{equation*}
    \mathsf{T}_q\Q\longrightarrow\Vspace_{v_q}\TQ,
    \qquad
    w\longmapsto w^{\Ver}_{v_q},
\end{equation*}
is a linear isomorphism.

\subsubsection{Connection map}
The horizontal subspace is defined using the connector induced by $\LC$. Let
\begin{equation*}
    \xi\in \mathsf{T}_{v_q}\TQ.
\end{equation*}
Choose a smooth curve
\begin{equation*}
    s\longmapsto x(s)=V(s)\in\TQ
\end{equation*}
such that
\begin{equation*}
    x(0)=v_q,
    \qquad
    \dot x(0)=\xi,
\end{equation*}
where $V(s)\in \mathsf{T}_{q(s)}\Q$ is a vector field along the base curve $q(s)$.

\begin{Definition}[Connection map]
The connection map, or connector,
\begin{equation*}
    K_{v_q}:\mathsf{T}_{v_q}\TQ\longrightarrow \mathsf{T}_q\Q
\end{equation*}
is defined by
\begin{equation*}
    K_{v_q}(\xi)
    :=\left.\frac{\cD V}{\dd s}\right|_{s=0}
    =\left.\LC_{\dot q(s)}V(s)\right|_{s=0}.
\end{equation*}
This value depends only on the tangent vector $\xi$ and not on the representative curve $x(s)$; hence $K$ is a well-defined smooth bundle map.
\end{Definition}

For a vertical lift, the connector recovers the vector being lifted:
\begin{equation}
    K_{v_q}\!\left(w^{\Ver}_{v_q}\right)=w.
    \label{eq:connector-vertical}
\end{equation}
Indeed, for $s\mapsto(v_q+s w)$ the base point is constant and the covariant derivative reduces to ordinary differentiation in the fixed vector space $\mathsf{T}_q\Q$.

\subsubsection{Horizontal subspace and horizontal lift}

\begin{Definition}[Horizontal subspace]
The horizontal subspace determined by $\LC$ is
\begin{equation*}
    \Hspace_{v_q}\TQ
    :=\ker K_{v_q}.
\end{equation*}
Thus $\xi\in \mathsf{T}_{v_q}\TQ$ is horizontal precisely when the vector component $v_q$ is infinitesimally parallel transported along the corresponding base displacement.
\end{Definition}

To make this explicit, fix $w\in \mathsf{T}_q\Q$. Choose a smooth curve $q(s)$ with
\begin{equation*}
    q(0)=q,
    \qquad
    \dot q(0)=w,
\end{equation*}
and let $V(s)$ be the parallel transport of $v_q$ along $q(s)$, i.e.,
\begin{equation*}
    V(0)=v_q,
    \qquad
    \frac{\cD V}{\dd s}=0.
\end{equation*}
Then
\begin{equation*}
    w^{\Hor}_{v_q}
    :=\left.\frac{\dd}{\dd s}\right|_{s=0}V(s)
\end{equation*}
is horizontal and satisfies
\begin{equation}
    T\piQ\!\left(w^{\Hor}_{v_q}\right)=w,
    \qquad
    K\!\left(w^{\Hor}_{v_q}\right)=0.
    \label{eq:horizontal-lift-properties}
\end{equation}
At the level of first-order tangent data, $w^{\Hor}_{v_q}$ is independent of the chosen base curve realizing $w$.

\begin{Definition}[Horizontal lift]
The unique element $w^{\Hor}_{v_q}\in\Hspace_{v_q}\TQ$ satisfying
\begin{equation*}
    T\piQ\!\left(w^{\Hor}_{v_q}\right)=w
\end{equation*}
is called the horizontal lift of $w\in \mathsf{T}_q\Q$ at $v_q$.
\end{Definition}
The restriction
\begin{equation*}
    T\piQ\big|_{\Hspace_{v_q}\TQ}:
    \Hspace_{v_q}\TQ\longrightarrow \mathsf{T}_q\Q
\end{equation*}
is therefore a linear isomorphism.

\subsubsection{Direct-sum decomposition}

\begin{Proposition}[Horizontal--vertical decomposition]
For every $v_q\in\TQ$,
\begin{equation*}
    \mathsf{T}_{v_q}\TQ
    =\Hspace_{v_q}\TQ\oplus\Vspace_{v_q}\TQ.
\end{equation*}
Equivalently,
\begin{equation}
\begin{aligned}
    \Phi_{v_q}:
    \mathsf{T}_{v_q}\TQ&\longrightarrow \mathsf{T}_q\Q\oplus \mathsf{T}_q\Q,\\
    \xi&\longmapsto\bigl(T\piQ(\xi),K(\xi)\bigr)
\end{aligned}
    \label{eq:Phi-isomorphism}
\end{equation}
is a linear isomorphism.
\end{Proposition}

\begin{proof}
Given $a,b\in \mathsf{T}_q\Q$, define
\begin{equation*}
    \xi:=a^{\Hor}_{v_q}+b^{\Ver}_{v_q}.
\end{equation*}
Using \eqref{eq:horizontal-lift-properties} and \eqref{eq:connector-vertical},
\begin{equation*}
    T\piQ(\xi)=a,
    \qquad
    K(\xi)=b.
\end{equation*}
Hence $\Phi_{v_q}$ is surjective. If $\Phi_{v_q}(\xi)=(0,0)$, then the horizontal component has zero projection and is therefore zero, while the vertical component has zero connector and is therefore zero. Hence $\xi=0$, so $\Phi_{v_q}$ is injective. Its inverse is
\begin{equation*}
    \Phi_{v_q}^{-1}(a,b)
    =a^{\Hor}_{v_q}+b^{\Ver}_{v_q}.
\end{equation*}
\end{proof}

Consequently, every $\xi\in \mathsf{T}_{v_q}\TQ$ admits the unique decomposition
\begin{equation}
     {
    \xi
    =\bigl(T\piQ(\xi)\bigr)^{\Hor}
    +\bigl(K(\xi)\bigr)^{\Ver}.
    }
    \label{eq:xi-decomposition}
\end{equation}

\subsubsection{Tangent vector to a curve in
\texorpdfstring{$\mathsf{T}\mathbb S^2$}{TS²}}
Let
\begin{equation*}
    x(t)=v_{q(t)}\in\TQ
\end{equation*}
be a smooth curve, with $v_{q(t)}\in \mathsf{T}_{q(t)}\Q$. Applying \eqref{eq:xi-decomposition} to $\dot x(t)$ yields
\begin{equation}
     {
    \dot x(t)
    =\dot q(t)^{\Hor}
    +\left(\frac{\cD v_q}{\dd t}\right)^{\Ver}.
    }
    \label{eq:tangent-curve-splitting}
\end{equation}
Indeed,
\begin{equation*}
    T\piQ(\dot x)=\dot q,
    \qquad
    K(\dot x)=\frac{\cD v_q}{\dd t}.
\end{equation*}
Equation \eqref{eq:tangent-curve-splitting} is the intrinsic replacement for the Euclidean shorthand $\frac{\dd}{\dd t}(q,v)=(\dot q,\dot v)$.

\subsection{Deterministic Controlled Point-Mass Dynamics}

\subsubsection{Second-order intrinsic dynamics}
A deterministic point mass on $\Q=\Sph$ is described by a smooth configuration curve $q(t)\in\Q$ with velocity
\begin{equation*}
    v_{q(t)}:=\dot q(t)\in \mathsf{T}_{q(t)}\Q.
\end{equation*}
Its intrinsic acceleration is the covariant acceleration
\begin{equation*}
    \frac{\cD v_q}{\dd t}
    =\LC_{\dot q}v_q
    =\LC_{v_q}v_q.
\end{equation*}
We assume direct acceleration actuation, i.e., the control is an arbitrary tangent vector
\begin{equation*}
    u(t)\in \mathsf{T}_{q(t)}\Q,
\end{equation*}
and the controlled mechanical equation is
\begin{equation}
     {
    \LC_{v_q}v_q=u.
    }
    \label{eq:det-second-order}
\end{equation}
Together with $\dot q=v_q$, this is the intrinsic second-order point-mass model.

In the terminology of geometric mechanics, \eqref{eq:det-second-order} is a simple mechanical control system with zero potential and with the kinetic-energy metric equal to $\mathbb{G}$; see, e.g., \cite{BulloLewis2005}. It is fully actuated because the admissible control set at each $q$ spans the entire tangent space $\mathsf{T}_q\Q$. Equivalently, with any local frame $\{E_1,E_2\}$,
\begin{equation*}
    u=u^1E_1+u^2E_2,
\end{equation*}
and the two independent input directions span $\mathsf{T}_q\Sph$.

\subsubsection{First-order dynamics on
\texorpdfstring{$\mathsf{T}\mathbb S^2$}{TS²}}
Define
\begin{equation*}
    x(t)=v_{q(t)}\in\TQ.
\end{equation*}
Since $v_q=\dot q$ and $\cD v_q/\dd t=u$, \eqref{eq:tangent-curve-splitting} gives
\begin{equation}
 {
    \dot x=v_q^{\Hor}+u^{\Ver}.
    }
    \label{eq:det_sys_dynamics_split}
\end{equation}
Equivalently,
\begin{equation*}
    T\piQ(\dot x)=v_q,
    \qquad
    K(\dot x)=u.
\end{equation*}
The term $v_q^{\Hor}$ is the geodesic-spray component. Define
\begin{equation*}
    \mathcal S(v_q):=v_q^{\Hor}.
\end{equation*}
Then $\mathcal S\in\mathfrak X(\TQ)$ is the geodesic spray of $\LC$; its integral curves satisfy
\begin{equation*}
    \dot q=v_q,
    \qquad
    \LC_{v_q}v_q=0.
\end{equation*}

For a state-feedback control, it is convenient to regard $u$ as a smooth section of the pullback bundle,
\begin{equation*}
    u\in\Gamma(\piQ^*\mathsf{T}\Q),
    \qquad
    u(v_q)\in \mathsf{T}_q\Q.
\end{equation*}
Its vertical lift defines a genuine vector field $u^{\Ver}$ on $\TQ$, and the controlled vector field is
\begin{equation*}
    F_u(v_q)
    =\mathcal S(v_q)+u(v_q)^{\Ver}
    =v_q^{\Hor}+u^{\Ver}.
\end{equation*}
For a prescribed open-loop control $u(t)\in \mathsf{T}_{q(t)}\Q$, the same expression is interpreted as a time-dependent vector field along the state trajectory.

\begin{Remark}[Why $F(v_q)=(v_q,u)$ is only shorthand]
The pair $(v_q,u)$ belongs to $\mathsf{T}_q\Q\oplus \mathsf{T}_q\Q$, not literally to $\mathsf{T}_{v_q}\TQ$. The rigorous statement is
\begin{equation*}
    \Phi_{v_q}\bigl(F_u(v_q)\bigr)=(v_q,u),
\end{equation*}
where $\Phi_{v_q}=(T\piQ,K)$ is the connection-dependent isomorphism in \eqref{eq:Phi-isomorphism}. Thus $F_u=v_q^{\Hor}+u^{\Ver}$ is the intrinsic realization of the informal ordered pair $(v_q,u)$.
\end{Remark}

\section{Problem Formulation}
\label{sec:Problem_formulation}
The objective of this paper is twofold. First, starting from the fully actuated deterministic point-mass dynamics~\eqref{eq:det_sys_dynamics_split}, we seek an entirely intrinsic stochastic state model on $\mathsf{T}\Sph$ in which random perturbations of configuration and velocity are represented by the horizontal and vertical geometry of the tangent bundle. Second, for a smooth Markov state-feedback control, we seek a Fokker-Planck PDE governing the probability density of this stochastic system with respect to the Riemannian volume measure induced by the Sasaki metric on the tangent bundle. Both constructions must be global geometric objects, even though local orthonormal frames are used to write representative SDEs.

\begin{Problem}
Let $(\Sph,\Grnd)$ be the unit two-sphere equipped with its standard round Riemannian metric $\Grnd$ and associated Levi-Civita connection $\LC$. With $\mathsf{T}\Sph$ being the state manifold, let $x=v_q \in \mathsf{T}\Sph$ denote the state of the system, and let $u\in\Gamma(\piQ^*\mathsf{T}\Sph)$ be a smooth fully actuating state-feedback acceleration field. Consider the deterministic system dynamics~\eqref{eq:det_sys_dynamics_split} of the point-mass system on $(\Sph,\Grnd)$. Let $\sigma_q\ge0$ denote the constant noise intensity in horizontal, i.e., position, directions and $\sigma_v\ge0$ the constant noise intensity in vertical, i.e., velocity, directions. Construct and characterize an intrinsic geometric model
for the controlled stochastic dynamics of the point-mass system and the corresponding Fokker-Planck equation governing the  probability density evolution of the system. 
\end{Problem}

More precisely, we seek a stochastic model that recovers the original deterministic dynamics in the absence of noise and incorporates independent isotropic perturbations along the horizontal and vertical subspaces of the tangent bundle. 
For a prescribed  initial probability density $\rho_0$ on $\mathsf{T}\Sph$, we aim to derive the associated Fokker-Planck equation by taking the formal adjoint of the infinitesimal generator associated with the SDE of the stochastic system with respect to the Sasaki Riemannian volume measure, expressing the resulting equation in an intrinsic double-divergence form. Furthermore, we require that the stochastic model and its corresponding Fokker-Planck equation be independent of the choice of local orthonormal frames 
and that the associated probability density remain normalized for all time under suitable regularity and decay assumptions.

\section{Controlled SDE for the Stochastic Dynamics of a Point Mass on
\texorpdfstring{$\mathbb{S}^2$}{S²}}
\label{sec:controlled_sde}

\subsection{The Sasaki Metric on
\texorpdfstring{$\mathsf{T}\mathbb S^2$}{TS²}}
The stochastic state equation evolves on $\TQ$, so a Riemannian metric on the state manifold is required in order to formulate Brownian motion intrinsically. The canonical metric associated with $(\Q,\mathbb{G},\LC)$ is the Sasaki metric.

\begin{Definition}[Sasaki metric]
For $\xi,\eta\in \mathsf{T}_{v_q}\TQ$, define
\begin{equation*}
     {
    \Sasaki_{v_q}(\xi,\eta)
    :=\mathbb{G}_q\bigl(T\piQ(\xi),T\piQ(\eta)\bigr)
    +\mathbb{G}_q\bigl(K(\xi),K(\eta)\bigr).
    }
\end{equation*}
\end{Definition}
By the horizontal--vertical decomposition,
\begin{align*}
    \Sasaki(a^{\Hor},b^{\Hor})&=\mathbb{G}(a,b),
    \\ 
    \Sasaki(a^{\Ver},b^{\Ver})&=\mathbb{G}(a,b),
    \\ 
    \Sasaki(a^{\Hor},b^{\Ver})&=0.
\end{align*}
Hence the horizontal and vertical distributions are orthogonal and each is isometric to $\mathsf{T}_q\Q$. Let $\LCs$ denote the Levi--Civita connection of $\Sasaki$ on $\TQ$.

\subsubsection{Lifted orthonormal frame}
Let $U\subset\Q$ be an open set on which a smooth local $\mathbb{G}$-orthonormal frame
\begin{equation*}
    \{E_1,E_2\},
    \qquad
    \mathbb{G}(E_i,E_j)=\delta_{ij}
\end{equation*}
is defined. On $\piQ^{-1}(U)\subset\TQ$, its horizontal and vertical lifts form
\begin{equation}
     {
    \mathcal E
    :=\{E_1^{\Hor},E_2^{\Hor},E_1^{\Ver},E_2^{\Ver}\},
    }
    \label{eq:sasaki-ON-frame}
\end{equation}
which is a local $\Sasaki$-orthonormal frame on the four-dimensional state manifold $\TQ$.

\subsubsection{Levi--Civita connection of the Sasaki metric}
For vector fields $X,Y\in\mathfrak X(\Q)$, the Levi--Civita connection $\LCs$ of the Sasaki metric satisfies the standard identities \cite{Sasaki1958,YanoIshihara1973}
\begin{align*}
    \LCs_{X^{\Hor}}Y^{\Hor}
    &=(\LC_XY)^{\Hor}
    -\frac12\bigl(\Riem(X,Y)v_q\bigr)^{\Ver},
    \\ 
    \LCs_{X^{\Hor}}Y^{\Ver}
    &=(\LC_XY)^{\Ver}
    +\frac12\bigl(\Riem(v_q,Y)X\bigr)^{\Hor},
    \\ 
    \LCs_{X^{\Ver}}Y^{\Hor}
    &=\frac12\bigl(\Riem(v_q,X)Y\bigr)^{\Hor},
    \\ 
    \LCs_{X^{\Ver}}Y^{\Ver}
    &=0.
\end{align*}
Setting $X=Y=E_i$ and using $\Riem(E_i,E_i)=0$ yields
\begin{equation}
     {
    \LCs_{E_i^{\Hor}}E_i^{\Hor}
    =(\LC_{E_i}E_i)^{\Hor},
    }
    \label{eq:sasaki-diagonal-H}
\end{equation}
and
\begin{equation}
     {
    \LCs_{E_i^{\Ver}}E_i^{\Ver}=0.
    }
    \label{eq:sasaki-diagonal-V}
\end{equation}
These two identities are the key simplifications used below.

\subsection{Intrinsic Brownian Motion on a Riemannian Manifold}
The stochastic construction is based on~\cite[Theorem~4]{LeeChirikjian2025}. Let $(\mathsf{M},h)$ be an $n$-dimensional Riemannian manifold with Levi--Civita connection $\nabla^{\mathsf{M}}$, and let $\{Z_1,\ldots,Z_n\}$ be a local $h$-orthonormal frame. Their intrinsic Stratonovich equation for Brownian motion is
\begin{equation}
    \dd Y_t
    =-\frac12\sum_{a=1}^{n}\nabla^{\mathsf{M}}_{Z_a}Z_a(Y_t)\,\dd t
    +\sum_{a=1}^{n}Z_a(Y_t)\circ\dd W_t^a.
    \label{eq:theorem4-general}
\end{equation}
The corresponding infinitesimal generator is
\begin{equation*}
    \frac12\Delta_{\mathsf{M}},
\end{equation*}
where $\Delta_{\mathsf{M}}$ denotes the Laplace--Beltrami operator of $h$. Thus the Stratonovich drift
\begin{equation*}
    -\frac12\sum_a\nabla^{\mathsf{M}}_{Z_a}Z_a
\end{equation*}
is precisely the geometric correction ensuring that the second-order generator is $(1/2)\Delta_{\mathsf{M}}$.

If all noise fields are multiplied by a constant $\sigma$, then
\begin{equation}
    -\frac12\sum_a\nabla^{\mathsf{M}}_{\sigma Z_a}(\sigma Z_a)
    =-\frac{\sigma^2}{2}\sum_a\nabla^{\mathsf{M}}_{Z_a}Z_a.
    \label{eq:sigma-correction}
\end{equation}
The quadratic scaling in \eqref{eq:sigma-correction} uses the assumption that $\sigma$ is constant. State-dependent noise intensities would contribute additional derivatives of the intensity fields.

\begin{Remark}[Role of the Brownian-motion result]
The intrinsic Brownian-motion construction \eqref{eq:theorem4-general} and its interpretation are taken from \cite[Theorem~4]{LeeChirikjian2025}; that work also gives an explicit construction on $\mathbb S^2$. The horizontal--vertical decomposition and Sasaki geometry used here are additional tangent-bundle structures required because the present state manifold is $\mathsf{T}\mathbb S^2$, not merely $\mathbb S^2$.
\end{Remark}

\subsection{Probabilistic Setup and Control}
Let
\begin{equation*}
    (\Omega,\mathcal F,(\mathcal F_t)_{t\ge0},\mathbb P)
\end{equation*}
be a filtered probability space satisfying the usual conditions. Let
\begin{equation*}
    W^{q,1},\ W^{q,2},\ W^{v,1},\ W^{v,2}
\end{equation*}
be four independent standard real-valued Wiener processes adapted to $(\mathcal F_t)$. Let
\begin{equation*}
    X_t=v_{q_t}\in\TQ
\end{equation*}
be the stochastic state. We introduce two nonnegative constant noise intensities,
\begin{equation*}
    \sigma_q\ge0,
    \qquad
    \sigma_v\ge0,
\end{equation*}
where $\sigma_q$ determines the intensity in horizontal, i.e., position, directions and $\sigma_v$ determines the intensity in vertical, i.e., velocity, directions.

The control satisfies
\begin{equation*}
    u_t\in \mathsf{T}_{q_t}\Q.
\end{equation*}
For a Markov state-feedback model we may take a smooth section
\begin{equation*}
    u:\TQ\longrightarrow \mathsf{T}\Q,
    \qquad
    u(v_q)\in \mathsf{T}_q\Q,
\end{equation*}
of $\piQ^*\mathsf{T}\Q$, whose vertical lift is denoted by $u^{\Ver}(v_q)$.

\subsection{Geometric Meaning of Position and Velocity Noise}

\subsubsection{Horizontal noise as position noise}
For $w^{\Hor}\in\Hspace_{v_q}\TQ$,
\begin{equation*}
    T\piQ(w^{\Hor})=w,
    \qquad
    K(w^{\Hor})=0.
\end{equation*}
Hence a horizontal displacement moves the base point $q$ in direction $w$ while parallel transporting the velocity vector $v_q$ to first order. It is therefore the intrinsic mechanism for introducing random displacement of the configuration without imposing an independent random fiber displacement. The horizontal position-noise fields are
\begin{equation*}
    \mathcal H_i:=E_i^{\Hor},
    \qquad i=1,2.
\end{equation*}

\subsubsection{Vertical noise as velocity noise}
For $w^{\Ver}\in\Vspace_{v_q}\TQ$,
\begin{equation*}
    T\piQ(w^{\Ver})=0,
    \qquad
    K(w^{\Ver})=w.
\end{equation*}
Thus a vertical displacement leaves $q$ fixed to first order and changes $v_q$ within the fiber $\mathsf{T}_q\Q$. The velocity-noise fields are
\begin{equation*}
    \mathcal V_i:=E_i^{\Ver},
    \qquad i=1,2.
\end{equation*}

\subsection{Controlled Stratonovich SDE on the State Manifold}
We now combine the deterministic controlled drift $v_q^{\Hor}+u^{\Ver}$ with horizontal and vertical diffusion fields, using the geometric correction principle of \eqref{eq:theorem4-general}. The diffusion vector fields are
\begin{equation*}
    \sigma_qE_1^{\Hor},\quad
    \sigma_qE_2^{\Hor},\quad
    \sigma_vE_1^{\Ver},\quad
    \sigma_vE_2^{\Ver}.
\end{equation*}
The associated Stratonovich correction on $(\TQ,\Sasaki)$ is therefore
\begin{equation*}
    -\frac{\sigma_q^2}{2}\sum_{i=1}^2
    \LCs_{E_i^{\Hor}}E_i^{\Hor}
    -\frac{\sigma_v^2}{2}\sum_{i=1}^2
    \LCs_{E_i^{\Ver}}E_i^{\Ver}.
\end{equation*}
Accordingly, the controlled stochastic state equation is
\begin{equation}
 {
\begin{aligned}
    \dd X_t
    &={}\Bigg[
        v_{q_t}^{\Hor}+u_t^{\Ver}
        -\frac{\sigma_q^2}{2}\sum_{i=1}^2
          \LCs_{E_i^{\Hor}}E_i^{\Hor}
        -\frac{\sigma_v^2}{2}\sum_{i=1}^2
          \LCs_{E_i^{\Ver}}E_i^{\Ver}
    \Bigg]\dd t
    \\
    &+\sigma_q\sum_{i=1}^2
        E_i^{\Hor}(X_t)\circ\dd W_t^{q,i}
    +\sigma_v\sum_{i=1}^2
        E_i^{\Ver}(X_t)\circ\dd W_t^{v,i}.
\end{aligned}}
\label{eq:stochastic-state-full}
\end{equation}
Every drift and diffusion term in \eqref{eq:stochastic-state-full} belongs to $\mathsf{T}_{X_t}\mathsf{T}\Sph$; hence the equation is formulated directly on the state manifold.

\subsubsection{Simplification using the Sasaki connection}
Using \eqref{eq:sasaki-diagonal-H}--\eqref{eq:sasaki-diagonal-V},
\begin{equation*}
    \LCs_{E_i^{\Hor}}E_i^{\Hor}
    =(\LC_{E_i}E_i)^{\Hor},
    \qquad
    \LCs_{E_i^{\Ver}}E_i^{\Ver}=0.
\end{equation*}
Therefore \eqref{eq:stochastic-state-full} becomes
\begin{equation}
 {
\begin{aligned}
    \dd X_t
    &={}\left[
        v_{q_t}^{\Hor}+u_t^{\Ver}
        -\frac{\sigma_q^2}{2}\sum_{i=1}^2
        (\LC_{E_i}E_i)^{\Hor}
    \right]\dd t
    \\
    &+\sigma_q\sum_{i=1}^2
        E_i^{\Hor}(X_t)\circ\dd W_t^{q,i}
    +\sigma_v\sum_{i=1}^2
        E_i^{\Ver}(X_t)\circ\dd W_t^{v,i}.
\end{aligned}}
\label{eq:stochastic-state-simplified}
\end{equation}
Equation \eqref{eq:stochastic-state-simplified} is the principal intrinsic controlled stochastic model of this paper. The horizontal correction survives because the horizontal lifts of a local orthonormal frame are generally not parallel with respect to $\LCs$, whereas the direct vertical correction vanishes because $\LCs_{X^{\Ver}}Y^{\Ver}=0$.

\subsection{Separated Position and Covariant-Velocity Equations}
Equation \eqref{eq:stochastic-state-simplified} can be decomposed into base and fiber components using $T\piQ$ and $K$. Because the equation is in Stratonovich form, smooth maps obey the ordinary chain rule.

Using
\begin{equation*}
\begin{gathered}
    T\piQ(v_q^{\Hor})=v_q,
    \qquad T\piQ(u^{\Ver})=0,\\
    T\piQ(E_i^{\Hor})=E_i,
    \qquad T\piQ(E_i^{\Ver})=0,\\
    T\piQ\bigl((\LC_{E_i}E_i)^{\Hor}\bigr)=\LC_{E_i}E_i,
\end{gathered}
\end{equation*}
we obtain the position equation

\begin{equation*}
\begin{aligned}
    \dd q_t
    ={}&\left[
        v_{q_t}-\frac{\sigma_q^2}{2}\sum_{i=1}^2
        \LC_{E_i}E_i(q_t)
    \right]\dd t\\
    &+\sigma_q\sum_{i=1}^2E_i(q_t)\circ\dd W_t^{q,i}.
\end{aligned}
\end{equation*}

For the velocity component, define the covariant Stratonovich differential of the tangent-vector process $v_{q_t}$ along $q_t$ by
\begin{equation*}
     {
    \cD v_{q_t}:=K_{X_t}(\circ\dd X_t).
    }
\end{equation*}
This is the stochastic counterpart of the deterministic identity $K(\dot x)=\cD v_q/\dd t$. Since
\begin{equation*}
\begin{gathered}
    K(v_q^{\Hor})=0,
    \qquad K(u^{\Ver})=u,\\
    K(E_i^{\Hor})=0,
    \qquad K(E_i^{\Ver})=E_i,\\
    K\bigl((\LC_{E_i}E_i)^{\Hor}\bigr)=0,
\end{gathered}
\end{equation*}
we obtain
\begin{equation}
 {
    \cD v_{q_t}
    =u_t\,\dd t
    +\sigma_v\sum_{i=1}^2E_i(q_t)\circ\dd W_t^{v,i}.
}
\label{eq:velocity-SDE}
\end{equation}
Thus the stochastic system is equivalently
\begin{equation}
 {
\begin{aligned}
    \dd q_t
    &={}
    \left[
        v_{q_t}-\frac{\sigma_q^2}{2}\sum_{i=1}^2
        \LC_{E_i}E_i(q_t)
    \right]\dd t\\
    &\quad
    +\sigma_q\sum_{i=1}^2E_i(q_t)\circ\dd W_t^{q,i},\\
    \cD v_{q_t}
    &={}
    u_t\,\dd t
    +\sigma_v\sum_{i=1}^2E_i(q_t)\circ\dd W_t^{v,i}.
\end{aligned}}
\label{eq:split-stochastic-system}
\end{equation}

\subsubsection{Interpretation of the horizontal noise}
No explicit $\sigma_q$ term appears on the right-hand side of \eqref{eq:velocity-SDE}. This does not mean that the random motion of $q_t$ is unrelated to $v_{q_t}$. Rather, the horizontal lift is defined precisely by the condition
\begin{equation*}
    K(E_i^{\Hor})=0.
\end{equation*}
Thus, when $q_t$ is displaced by a horizontal stochastic increment, $v_{q_t}$ is simultaneously moved by stochastic parallel transport. In a coordinate representation this transport would appear through connection-dependent terms in the components of $v_{q_t}$; in the intrinsic formulation those terms are absorbed into the covariant differential $\cD v_{q_t}$.

The two classes of noise therefore have distinct geometric roles: $E_i^{\Hor}\circ\dd W^{q,i}$ randomly displaces position while parallel transporting the velocity state, whereas $E_i^{\Ver}\circ\dd W^{v,i}$ randomly changes the velocity inside the current tangent space $\mathsf{T}_{q_t}\Q$.

\subsubsection{Meaning of
\texorpdfstring{$v_{q_t}$}{vqt} when
\texorpdfstring{$\sigma_q>0$}{sigmaq > 0}}
When $\sigma_q=0$, the position process has finite variation and
\begin{equation*}
    \dd q_t=v_{q_t}\,\dd t,
\end{equation*}
so $v_{q_t}$ retains the classical mechanical meaning $v_{q_t}=\dot q_t$ wherever the derivative is defined. When $\sigma_q>0$, however, $q_t$ contains a Brownian component and its sample paths are almost surely nowhere classically differentiable. Consequently,
\begin{equation*}
    v_{q_t}=\dot q_t
\end{equation*}
cannot hold pathwise. In \eqref{eq:split-stochastic-system}, $v_{q_t}$ is instead the tangent drift-velocity state: it determines the deterministic first-order transport in the finite-variation part of the configuration dynamics, while the configuration also undergoes independent stochastic displacement. This distinction is unavoidable when idealized white noise acts directly in position.

\subsection{Infinitesimal Generator}
Assume that $u=u(v_q)$ is a smooth state feedback, so the diffusion is Markovian. Let $f\in C^\infty(\TQ)$. For a Stratonovich SDE
\begin{equation*}
    \dd X=A(X)\,\dd t+\sum_a B_a(X)\circ\dd W^a,
\end{equation*}
the infinitesimal generator is
\begin{equation}
    \mathcal A f
    =A[f]+\frac12\sum_a B_a[B_a[f]].
    \label{eq:general-strat-generator}
\end{equation}
Applying \eqref{eq:general-strat-generator} to \eqref{eq:stochastic-state-full}, and using the Hessian $\Hess^{\mathrm S}$ associated with $\LCs$, yields
\begin{equation}
 {
\begin{aligned}
    \mathcal A^u f
    ={}&(v_q^{\Hor}+u^{\Ver})[f]
    +\frac{\sigma_q^2}{2}\sum_{i=1}^2
      \Hess^{\mathrm S}f(E_i^{\Hor},E_i^{\Hor})\\
    &+\frac{\sigma_v^2}{2}\sum_{i=1}^2
      \Hess^{\mathrm S}f(E_i^{\Ver},E_i^{\Ver}).
\end{aligned}}
\label{eq:generator-HV-Hessian}
\end{equation}
Define
\begin{align*}
    \Delta_{\Hor}f
    &:=\sum_{i=1}^2
      \Hess^{\mathrm S}f(E_i^{\Hor},E_i^{\Hor}),
    \\ 
    \Delta_{\Ver}f
    &:=\sum_{i=1}^2
      \Hess^{\mathrm S}f(E_i^{\Ver},E_i^{\Ver}).
\end{align*}
These are traces over the horizontal and vertical subspaces with respect to the Sasaki metric, and hence do not depend on the orthonormal frame used to evaluate the traces. Therefore
\begin{equation*}
     {
    \mathcal A^u
    =v_q^{\Hor}+u^{\Ver}
    +\frac{\sigma_q^2}{2}\Delta_{\Hor}
    +\frac{\sigma_v^2}{2}\Delta_{\Ver},
    }
\end{equation*}
where the first two terms are understood as first-order differential operators.

\subsubsection{Equal position and velocity noise intensities}
If
\begin{equation*}
    \sigma_q=\sigma_v=\sigma,
\end{equation*}
then \eqref{eq:sasaki-ON-frame} is a complete local Sasaki-orthonormal frame, and
\begin{equation*}
    \Delta^{\mathrm S}=\Delta_{\Hor}+\Delta_{\Ver},
\end{equation*}
where $\Delta^{\mathrm S}$ is the Laplace--Beltrami operator of $(\TQ,\Sasaki)$. Consequently,
\begin{equation}
     {
    \mathcal A^u f
    =(v_q^{\Hor}+u^{\Ver})[f]
    +\frac{\sigma^2}{2}\Delta^{\mathrm S}f.
    }
    \label{eq:generator-equal-noise}
\end{equation}
Thus, for equal noise intensities, the stochastic perturbation is an isotropic Brownian perturbation on the full Sasaki state manifold, superimposed on the deterministic controlled drift. If $\sigma_q\neq\sigma_v$, the model remains intrinsic and isotropic within each of the horizontal and vertical subbundles, but it is anisotropic with respect to the unweighted Sasaki metric on the full state space.

\subsection{Intrinsic \Ito{} Representation on the State Manifold}
Although the Stratonovich representation is the geometrically natural one, it is useful to record the corresponding intrinsic \Ito{} form relative to $\LCs$. For the diffusion fields in \eqref{eq:stochastic-state-full}, the \Ito{}--Stratonovich correction is
\begin{equation*}
    \frac{\sigma_q^2}{2}\sum_{i=1}^2
    \LCs_{E_i^{\Hor}}E_i^{\Hor}
    +\frac{\sigma_v^2}{2}\sum_{i=1}^2
    \LCs_{E_i^{\Ver}}E_i^{\Ver}.
\end{equation*}
This cancels the geometric correction already included in the Stratonovich drift. Hence the intrinsic \Ito{} equation on $(\TQ,\Sasaki)$ is
\begin{equation}
 {
\begin{aligned}
    \dd X_t
    ={}&(v_{q_t}^{\Hor}+u_t^{\Ver})\,\dd t
    +\sigma_q\sum_{i=1}^2E_i^{\Hor}(X_t)\,\dd W_t^{q,i}\\
    &+\sigma_v\sum_{i=1}^2E_i^{\Ver}(X_t)\,\dd W_t^{v,i}.
\end{aligned}}
\label{eq:intrinsic-Ito-state}
\end{equation}
Equation \eqref{eq:intrinsic-Ito-state} is an intrinsic \Ito{} representation defined relative to the Sasaki connection; it is not an ambient Euclidean \Ito{} equation.

\subsection{Local Frames and Global Meaning}
The sphere $\Sph$ does not admit a global smooth nonvanishing orthonormal frame. Thus no single smooth pair $\{E_1,E_2\}$ exists on all of $\Sph$, and the displayed Stratonovich equations are local frame representations on open sets $U\subset\Sph$.

This does not make the diffusion law frame dependent. If $\{F_1,F_2\}$ is another local orthonormal frame, then
\begin{equation*}
    F_i=\sum_{j=1}^2O_{ij}E_j,
\end{equation*}
for an orthogonal matrix-valued map $O$. Since trace is invariant under orthogonal changes of basis,
\begin{align*}
    \sum_i\Hess^{\mathrm S}f(F_i^{\Hor},F_i^{\Hor})
    &=\sum_i\Hess^{\mathrm S}f(E_i^{\Hor},E_i^{\Hor}),\\
    \sum_i\Hess^{\mathrm S}f(F_i^{\Ver},F_i^{\Ver})
    &=\sum_i\Hess^{\mathrm S}f(E_i^{\Ver},E_i^{\Ver}).
\end{align*}
Therefore the generator is globally defined. This is the same frame-independence principle emphasized for intrinsic Brownian motion in \cite{LeeChirikjian2025}: local Stratonovich representations may change with the orthonormal frame, but the induced second-order diffusion operator is intrinsic.

\subsection{Consistency Checks and Limiting Cases}

\subsubsection{No stochastic noise}
If
\begin{equation*}
    \sigma_q=\sigma_v=0,
\end{equation*}
then \eqref{eq:stochastic-state-simplified} reduces to
\begin{equation*}
    \dot X=v_q^{\Hor}+u^{\Ver},
\end{equation*}
and the split equations become
\begin{equation*}
    \dot q=v_q,
    \qquad
    \frac{\cD v_q}{\dd t}=u.
\end{equation*}
Thus the deterministic point-mass model is recovered exactly.

\subsubsection{Velocity noise only}
If $\sigma_q=0$ and $\sigma_v>0$, then
\begin{align*}
    \dd q_t&=v_{q_t}\,\dd t,\\
    \cD v_{q_t}&=u_t\,\dd t
    +\sigma_v\sum_{i=1}^2E_i(q_t)\circ\dd W_t^{v,i}.
\end{align*}
The position remains of finite variation, while the acceleration is stochastically perturbed in every tangent direction.

\subsubsection{Position noise only}
If $\sigma_q>0$ and $\sigma_v=0$, then
\begin{align*}
    \dd q_t
    &=\left[
    v_{q_t}-\frac{\sigma_q^2}{2}\sum_i\LC_{E_i}E_i(q_t)
    \right]\dd t\\ \notag
    &\quad\quad+\sigma_q\sum_iE_i(q_t)\circ\dd W_t^{q,i},\\ 
    \cD v_{q_t}&=u_t\,\dd t.
\end{align*}
Thus the position receives intrinsic Brownian displacement while the velocity is covariantly transported except for its deterministic control-induced change.

\subsubsection{Equal isotropic state noise}
If $\sigma_q=\sigma_v=\sigma$, then
\begin{equation*}
    \mathcal A^u
    =v_q^{\Hor}+u^{\Ver}+\frac{\sigma^2}{2}\Delta^{\mathrm S}.
\end{equation*}
This is the direct controlled analogue on $(\mathsf{T}\Sph,\Sasaki)$ of the intrinsic Brownian-motion construction of \cite{LeeChirikjian2025}.

\subsection{Collected Final Form of the Model}
For reference, the deterministic and stochastic equations are collected here.

\subsubsection{Deterministic point mass}
The state is
\begin{equation*}
    x=v_q\in \mathsf{T}\Sph,
\end{equation*}
with
\begin{equation*}
    v_q=\dot q,
    \qquad
    u=\LC_{v_q}v_q.
\end{equation*}
The first-order intrinsic dynamics are
\begin{equation*}
     {
    \dot x=v_q^{\Hor}+u^{\Ver}.
    }
\end{equation*}

\subsubsection{Stochastic point mass with position and velocity noise}
Let $\{E_1,E_2\}$ be any local orthonormal frame on $\Sph$. With the Sasaki metric on $\mathsf{T}\Sph$ and constant $\sigma_q,\sigma_v\ge0$, the intrinsic Stratonovich SDE is
\begin{equation*}
 {
\begin{aligned}
    \dd X_t
    &={}\left[
        v_{q_t}^{\Hor}+u_t^{\Ver}
        -\frac{\sigma_q^2}{2}\sum_{i=1}^2
        (\LC_{E_i}E_i)^{\Hor}
    \right]\dd t\\
    &+\sigma_q\sum_{i=1}^2E_i^{\Hor}(X_t)\circ\dd W_t^{q,i}
    +\sigma_v\sum_{i=1}^2E_i^{\Ver}(X_t)\circ\dd W_t^{v,i}.
\end{aligned}}
\end{equation*}
Equivalently,
\begin{equation*}
 {
\begin{aligned}
    \dd q_t
    &={}
    \left[
        v_{q_t}-\frac{\sigma_q^2}{2}\sum_{i=1}^2
        \LC_{E_i}E_i(q_t)
    \right]\dd t \\ 
    &\quad\quad+\sigma_q\sum_{i=1}^2E_i(q_t)\circ\dd W_t^{q,i},\\
    \cD v_{q_t}
    &={}
    u_t\,\dd t
    +\sigma_v\sum_{i=1}^2E_i(q_t)\circ\dd W_t^{v,i}.
\end{aligned}}
\end{equation*}
The covariant Stratonovich differential is
\begin{equation*}
    \cD v_{q_t}=K_{X_t}(\circ\dd X_t).
\end{equation*}
For equal intensities $\sigma_q=\sigma_v=\sigma$, the generator is
\begin{equation*}
     {
    \mathcal A^u f
    =(v_q^{\Hor}+u^{\Ver})[f]
    +\frac{\sigma^2}{2}\Delta^{\mathrm S}f,
    \qquad f\in C^\infty(\mathsf{T}\Sph).
    }
\end{equation*}

\section{Intrinsic Fokker--Planck Equation on
\texorpdfstring{$\mathsf{T}\mathbb S^2$}{TS²}}
\label{sec:intrinsic-FP}

We now derive the deterministic partial differential equation governing the probability density of the stochastic state process $X_t=v_{q_t}\in \mathsf{T}\Sph$.  The derivation is carried out directly on the four-dimensional Riemannian state manifold
\[
    \mathsf{M}:=\mathsf{T}\Q=\mathsf{T}\Sph
\]
equipped with the Sasaki metric $\Sasaki$ and its Levi--Civita connection $\LCs$.  The treatment follows the global tensor-field formulation of the Fokker--Planck equation developed in~\cite{LeeChirikjianFP2026}.  In particular, we use their intrinsic adjoint formulation for Stratonovich SDEs and their diffusion-tensor formulation for intrinsic It\^o SDEs.

\subsection{The Sasaki volume measure and probability density}

Because $(\mathsf{M},\Sasaki)$ is a Riemannian manifold, it carries the Riemannian volume measure induced by $\Sasaki$.  We denote this measure by
\[
    \muS.
\]
When an orientation is fixed, the same notation may be used for the corresponding Riemannian volume form.  All probability densities below are densities \emph{with respect to this intrinsic Sasaki volume measure}; no coordinate Lebesgue measure and no ambient Euclidean measure are used.

Suppose that the law of $X_t$ is absolutely continuous with respect to $\muS$.  Its probability density is denoted by
\[
    \rho_t(x)=\rho(t,x),
    \qquad
    x=v_q\in \mathsf{M}.
\]
Thus, for every Borel set $A\subset \mathsf{M}$,
\begin{equation*}
    \mathbb P(X_t\in A)
    =
    \int_A \rho_t(x)\,\dd\muS(x).
\end{equation*}
A probability density satisfies
\begin{equation*}
    \rho_t(x)\ge0,
    \qquad
    \int_{\mathsf{M}} \rho_t\,\dd\muS=1.
\end{equation*}

The state manifold $\mathsf{T}\Sph$ is noncompact because every tangent fiber is isomorphic to $\mathbb R^2$.  In the formal derivations below we therefore use compactly supported smooth test functions.  When integrations by parts are applied directly to $\rho_t$, it is understood that $\rho_t$ and the associated probability flux decay sufficiently rapidly in the fiber directions so that no boundary contribution appears ``at infinity.''  Equivalently, the Fokker--Planck equation can be interpreted in the usual weak sense against compactly supported test functions.  We also assume that the feedback $u=u(v_q)$ and the coefficients of the SDE satisfy the regularity and growth conditions needed for the state process to be well defined and nonexplosive on the time interval of interest.

\subsection{Intrinsic divergence on the Sasaki state manifold}

The divergence operator that appears in the forward equation is the divergence associated with the Sasaki metric.  For a smooth vector field
\[
    Y\in\mathfrak X(\mathsf{M}),
\]
define
\begin{equation}
     {
    \divS Y
    :=
    \trace\bigl(Z\mapsto\LCs_ZY\bigr).
    }
    \label{eq:sasaki-divergence-vector}
\end{equation}
If
\[
    \{\mathcal E_1,\ldots,\mathcal E_4\}
\]
is any local $\Sasaki$-orthonormal frame on $\mathsf{M}$, then
\begin{equation*}
    \divS Y
    =
    \sum_{a=1}^4
    \Sasaki\bigl(\LCs_{\mathcal E_a}Y,\mathcal E_a\bigr).
\end{equation*}
This expression is independent of the chosen orthonormal frame.

For scalar functions $r,s\in C_c^\infty(\mathsf{M})$, define the $L^2$ pairing
\begin{equation*}
    \langle r,s\rangle_{L^2(\muS)}
    :=
    \int_{\mathsf{M}} r(x)s(x)\,\dd\muS(x).
\end{equation*}
The intrinsic integration-by-parts identity used repeatedly below is
\begin{equation}
     {
    \langle r,Y[s]\rangle_{L^2(\muS)}
    =
    \left\langle-\divS(rY),s\right\rangle_{L^2(\muS)}.
    }
    \label{eq:IBP-vector}
\end{equation}
This is precisely the manifold identity used in~\cite[Lemma~1]{LeeChirikjianFP2026}.  For completeness, it follows from the product rule
\begin{equation*}
    \divS(rY)=Y[r]+r\,\divS Y
\end{equation*}
and the divergence theorem.  Since $r$ and $s$ are compactly supported,
\[
    \int_{\mathsf{M}} \divS(rsY)\,\dd\muS=0.
\]
Expanding the divergence of $rsY$ and rearranging gives \eqref{eq:IBP-vector}.

\subsection{Divergence of a contravariant tensor field}

The intrinsic It\^o Fokker--Planck equation will involve a symmetric tensor field of type $(2,0)$.  We therefore record the tensor divergence carefully.

Let
\[
    T\in\mathcal T^{(2,0)}(\mathsf{M})
\]
be a smooth contravariant two-tensor.  Let
\[
    \{\mathcal E_a\}_{a=1}^4
\]
be a local Sasaki-orthonormal frame and let
\[
    \{\Theta^a\}_{a=1}^4
\]
be its dual coframe.  The divergence of $T$ is the vector field $\divS T$ characterized by
\begin{equation*}
     {
    (\divS T)(\alpha)
    :=
    \sum_{a=1}^4
    \bigl(\LCs_{\mathcal E_a}T\bigr)
    \bigl(\alpha,\Theta^a\bigr),
    \qquad
    \alpha\in\Omega^1(\mathsf{M}).
    }
\end{equation*}
The right-hand side is independent of the chosen Sasaki-orthonormal frame.  This is the specialization to $(2,0)$ tensors of the tensor-divergence definition in~\cite{LeeChirikjianFP2026}.

The second divergence
\[
    \divS\bigl(\divS T\bigr)
\]
is then the ordinary divergence \eqref{eq:sasaki-divergence-vector} of the vector field $\divS T$.  Thus
\[
    \divS\bigl(\divS T\bigr)
\]
is a scalar field.

The tensor integration-by-parts identity of~\cite[Lemma~2]{LeeChirikjianFP2026}, applied twice, implies the following important relation.  For $f\in C_c^\infty(\mathsf{M})$ and a compactly supported symmetric $(2,0)$ tensor $T$,
\begin{equation}
     {
    \int_{\mathsf{M}}
    \bigl(T:\Hess^{\mathrm S}f\bigr)
    \,\dd\muS
    =
    \int_{\mathsf{M}}
    f\,
    \divS\bigl(\divS T\bigr)
    \,\dd\muS.
    }
    \label{eq:tensor-double-IBP}
\end{equation}
Here $\Hess^{\mathrm S}f$ is the $(0,2)$ Hessian associated with $\LCs$, and ``$:$'' denotes full contraction of the two tensor slots.

\subsection{Generator, adjoint, and the abstract forward equation}

We next recall why the Fokker--Planck equation is the adjoint evolution associated with the infinitesimal generator.  Let
\[
    \mathcal A^u:C^\infty(\mathsf{M})\longrightarrow C^\infty(\mathsf{M})
\]
denote the generator of the controlled Markov diffusion.  Its formal $L^2(\muS)$ adjoint $(\mathcal A^u)^*$ is defined by
\begin{equation*}
    \left\langle \rho,\mathcal A^u f\right\rangle_{L^2(\muS)}
    =
    \left\langle(\mathcal A^u)^*\rho,f\right\rangle_{L^2(\muS)}
\end{equation*}
for smooth compactly supported test functions.

Let $f\in C_c^\infty(\mathsf{M})$.  From the density representation,
\begin{equation*}
    \mathbb E[f(X_t)]
    =
    \int_{\mathsf{M}} f(x)\rho_t(x)\,\dd\muS(x).
\end{equation*}
Differentiating formally in time gives
\begin{equation*}
    \frac{\dd}{\dd t}\mathbb E[f(X_t)]
    =
    \left\langle\frac{\partial \rho_t}{\partial t},f\right\rangle_{L^2(\muS)}.
\end{equation*}
On the other hand, the defining property of the generator gives

\begin{equation*}
\begin{aligned}
    \frac{\dd}{\dd t}\mathbb E[f(X_t)]
    &=
    \mathbb E[(\mathcal A^u f)(X_t)]\\
    &=
    \langle \rho_t,\mathcal A^u f\rangle_{L^2(\muS)}\\
    &=
    \langle(\mathcal A^u)^*\rho_t,f\rangle_{L^2(\muS)}.
\end{aligned}
\end{equation*}

Since this holds for every $f\in C_c^\infty(\mathsf{M})$,
\begin{equation}
     {
    \frac{\partial \rho_t}{\partial t}
    =
    (\mathcal A^u)^*\rho_t.
    }
    \label{eq:abstract-FP}
\end{equation}
This is the abstract Fokker--Planck equation, exactly as developed in~\cite[Section~III-A]{LeeChirikjianFP2026}.

\subsection{Stratonovich Fokker--Planck equation for the controlled point mass}

We first apply the Stratonovich formulation because the stochastic mechanical model was constructed intrinsically in Stratonovich form.

Define the deterministic controlled mechanical drift
\begin{equation*}
    F_u
    :=
    v_q^\Hor+u^\Ver
    \in\mathfrak X(\mathsf{M}).
\end{equation*}
The actual Stratonovich drift in \eqref{eq:stochastic-state-simplified} is
\begin{equation*}
     {
    X_u^{\mathrm{Strat}}
    :=
    F_u
    -\frac{\sigma_q^2}{2}
    \sum_{i=1}^2
    (\LC_{E_i}E_i)^\Hor.
    }
\end{equation*}
Introduce the four diffusion vector fields
\begin{align*}
    B_i^q
    &:=
    \sigma_q E_i^\Hor,
    \qquad i=1,2,
    \\ 
    B_i^v
    &:=
    \sigma_v E_i^\Ver,
    \qquad i=1,2.
\end{align*}
Then the stochastic state equation can be written compactly as
\begin{equation}
\begin{aligned}
    \dd X_t
    ={}&
    X_u^{\mathrm{Strat}}(X_t)\,\dd t
    +
    \sum_{i=1}^2 B_i^q(X_t)\circ\dd W_t^{q,i}\\
    &+
    \sum_{i=1}^2 B_i^v(X_t)\circ\dd W_t^{v,i}.
\end{aligned}
\label{eq:compact-Strat-for-FP}
\end{equation}

For a Stratonovich SDE
\[
    \dd X=X(X)\,\dd t+\sum_\alpha B_\alpha(X)\circ\dd W^\alpha,
\]
~\cite[Theorem~2]{LeeChirikjianFP2026} gives the intrinsic forward equation
\[
    \frac{\partial \rho}{\partial t}
    =
    -\divS(\rho X)
    +
    \frac12
    \sum_\alpha
    \divS\!\left(
        \divS(\rho B_\alpha)\,B_\alpha
    \right).
\]
Applying this result to \eqref{eq:compact-Strat-for-FP} yields
\begin{equation}
 {
\begin{aligned}
    \frac{\partial \rho_t}{\partial t}
    ={}&
    -\divS\!\left(\rho_tX_u^{\mathrm{Strat}}\right)
    \\
    &+
    \frac12\sum_{i=1}^2
    \divS\!\left(
        \divS(\rho_tB_i^q)\,B_i^q
    \right)
    \\
    &+
    \frac12\sum_{i=1}^2
    \divS\!\left(
        \divS(\rho_tB_i^v)\,B_i^v
    \right).
\end{aligned}}
\label{eq:FP-Strat-compact}
\end{equation}
Because $\sigma_q$ and $\sigma_v$ are constants, this can equivalently be displayed as
\begin{equation}
 {
\begin{aligned}
    \frac{\partial \rho_t}{\partial t}
    ={}&
    -\divS\!\left[
        \rho_t\left(
        F_u
        -\frac{\sigma_q^2}{2}
        \sum_{i=1}^2(\LC_{E_i}E_i)^\Hor
        \right)
    \right]
    \\
    &+
    \frac{\sigma_q^2}{2}
    \sum_{i=1}^2
    \divS\!\left(
        \divS(\rho_tE_i^\Hor)\,E_i^\Hor
    \right)
    \\
    &+
    \frac{\sigma_v^2}{2}
    \sum_{i=1}^2
    \divS\!\left(
        \divS(\rho_tE_i^\Ver)\,E_i^\Ver
    \right).
\end{aligned}}
\label{eq:FP-Strat-expanded}
\end{equation}
Equation~\eqref{eq:FP-Strat-expanded} is an intrinsic Stratonovich Fokker--Planck equation on $(\mathsf{T}\Sph,\Sasaki)$.

\begin{remark}[Meaning of the nested divergence in the Stratonovich formula]
For a vector field $B$, the inner quantity
\[
    \divS(\rho B)
\]
is a scalar field.  Therefore
\[
    \divS(\rho B)\,B
\]
is again a vector field, and the outer divergence produces a scalar.  Thus every term in \eqref{eq:FP-Strat-compact} is intrinsically well defined on $\mathsf{M}$.
\end{remark}

\subsection{The intrinsic diffusion tensor on
\texorpdfstring{$\mathsf{T}\mathbb S^2$}{TS²}}

The tensor-field formulation of~\cite{LeeChirikjianFP2026} is especially useful for our model because the stochastic perturbation distinguishes the horizontal and vertical subbundles.

Define the horizontal and vertical contravariant metric tensors
\begin{align}
    \mathsf G_{\Hor}
    &:=
    \sum_{i=1}^2
    E_i^\Hor\otimes E_i^\Hor,
    \label{eq:GH-def}\\
    \mathsf G_{\Ver}
    &:=
    \sum_{i=1}^2
    E_i^\Ver\otimes E_i^\Ver.
    \label{eq:GV-def}
\end{align}
Although \eqref{eq:GH-def}--\eqref{eq:GV-def} are written using a local orthonormal frame on $\Q$, the tensors are globally defined.  Indeed, if
\[
    F_i=\sum_{j=1}^2O_{ij}E_j,
    \qquad
    O\in O(2),
\]
then
\[
    \sum_iF_i^\Hor\otimes F_i^\Hor
    =
    \sum_jE_j^\Hor\otimes E_j^\Hor,
\]
and the same computation holds for the vertical lifts.  Thus $\mathsf G_{\Hor}$ and $\mathsf G_{\Ver}$ depend only on the horizontal--vertical splitting and the Sasaki metric.

The diffusion tensor associated with \eqref{eq:intrinsic-Ito-state} is
\begin{equation}
\begin{aligned}
    \mathsf D
    &:=
    \sum_{i=1}^2B_i^q\otimes B_i^q
    +
    \sum_{i=1}^2B_i^v\otimes B_i^v
    \\
    &=
     {
    \sigma_q^2\mathsf G_{\Hor}
    +
    \sigma_v^2\mathsf G_{\Ver}.
    }
\end{aligned}
\label{eq:diffusion-tensor-HV}
\end{equation}
It is a smooth symmetric $(2,0)$ tensor field.  For every covector $\alpha\in \mathsf{T}_x^*\mathsf{M}$,
\begin{equation*}
\begin{aligned}
    \mathsf D(\alpha,\alpha)
    &=
    \sigma_q^2\sum_{i=1}^2
    \alpha(E_i^\Hor)^2
    +
    \sigma_v^2\sum_{i=1}^2
    \alpha(E_i^\Ver)^2
    \ge0,
\end{aligned}
\end{equation*}
so $\mathsf D$ is positive semidefinite.  It is precisely the intrinsic local covariance tensor of the stochastic perturbation, in the sense of~\cite[Theorem~4 and Remark~1]{LeeChirikjianFP2026}.

Because the lifted fields form a complete Sasaki-orthonormal frame,
\begin{equation}
     {
    \mathsf G_{\Hor}+\mathsf G_{\Ver}
    =
    (\Sasaki)^{-1},
    }
    \label{eq:cometric-splitting}
\end{equation}
where $(\Sasaki)^{-1}$ denotes the inverse Sasaki metric, viewed as a symmetric $(2,0)$ tensor.

\subsection{The generator in diffusion-tensor form}

The intrinsic It\^o representation \eqref{eq:intrinsic-Ito-state} has drift
\[
    \widetilde X_u=F_u=v_q^\Hor+u^\Ver
\]
and diffusion tensor $\mathsf D$ given by \eqref{eq:diffusion-tensor-HV}. The statement in~\cite[Theorem~4]{LeeChirikjianFP2026} states that the generator of an intrinsic It\^o SDE is
\[
    \mathcal A f
    =
    \widetilde X[f]
    +
    \frac12(\mathsf D:\Hess f).
\]
In our setting,
\begin{equation}
     {
    \mathcal A^u f
    =
    F_u[f]
    +
    \frac12
    \bigl(
        \mathsf D:\Hess^{\mathrm S}f
    \bigr).
    }
    \label{eq:generator-diffusion-tensor}
\end{equation}
Using \eqref{eq:diffusion-tensor-HV},
\begin{equation*}
\begin{aligned}
    \mathsf D:\Hess^{\mathrm S}f
    ={}&
    \sigma_q^2
    \sum_{i=1}^2
    \Hess^{\mathrm S}f(E_i^\Hor,E_i^\Hor)
    \\
    &+
    \sigma_v^2
    \sum_{i=1}^2
    \Hess^{\mathrm S}f(E_i^\Ver,E_i^\Ver).
\end{aligned}
\end{equation*}
Hence \eqref{eq:generator-diffusion-tensor} is exactly the generator already obtained in \eqref{eq:generator-HV-Hessian}.  The diffusion-tensor formulation therefore does not introduce a different stochastic process; it repackages the same second-order operator in a global tensorial form particularly well suited to taking the adjoint.

\subsection{Main intrinsic Fokker--Planck equation in double-divergence form}

We can now state the central forward equation.

\begin{theorem}[Intrinsic Fokker--Planck equation for the controlled stochastic point mass]
\label{thm:main-FP}
Let $X_t=v_{q_t}$ evolve according to the controlled stochastic dynamics \eqref{eq:stochastic-state-simplified}, with smooth Markov feedback $u=u(v_q)$ and constant noise intensities $\sigma_q,\sigma_v\ge0$.  Let $\rho_t$ be the probability density of $X_t$ with respect to the Sasaki volume measure $\muS$.  Under the regularity and decay assumptions stated above, $\rho_t$ satisfies
\begin{equation}
 {
    \frac{\partial \rho_t}{\partial t}
    =
    -\divS(\rho_tF_u)
    +
    \frac12
    \divS\!\left(
        \divS(\rho_t\mathsf D)
    \right),
}
\label{eq:main-FP-double-div}
\end{equation}
where
\begin{equation*}
    F_u=v_q^\Hor+u^\Ver
\end{equation*}
and
\begin{equation}
     {
    \mathsf D
    =
    \sigma_q^2\mathsf G_{\Hor}
    +
    \sigma_v^2\mathsf G_{\Ver}.
    }
    \label{eq:FP-D-repeat}
\end{equation}
Equivalently,
\begin{equation}
 {
\begin{aligned}
    \frac{\partial \rho_t}{\partial t}
    ={}&
    -\divS\!\left(
        \rho_t(v_q^\Hor+u^\Ver)
    \right)
    \\
    &+
    \frac{\sigma_q^2}{2}
    \divS\!\left(
        \divS(\rho_t\mathsf G_{\Hor})
    \right)
    \\
    &+
    \frac{\sigma_v^2}{2}
    \divS\!\left(
        \divS(\rho_t\mathsf G_{\Ver})
    \right).
\end{aligned}}
\label{eq:main-FP-HV}
\end{equation}
\end{theorem}

\begin{proof}
Let $f\in C_c^\infty(\mathsf{M})$.  From \eqref{eq:generator-diffusion-tensor},
\[
    \mathcal A^u f
    =
    F_u[f]
    +
    \frac12(\mathsf D:\Hess^{\mathrm S}f).
\]
Therefore
\begin{equation}
\begin{aligned}
    \langle \rho_t,\mathcal A^u f\rangle_{L^2(\muS)}
    ={}&
    \langle \rho_t,F_u[f]\rangle_{L^2(\muS)}
    \\
    &+
    \frac12
    \int_{\mathsf{M}}
    \rho_t\,
    \bigl(\mathsf D:\Hess^{\mathrm S}f\bigr)
    \,\dd\muS.
\end{aligned}
\label{eq:FP-proof-1}
\end{equation}
By the vector-field integration-by-parts identity \eqref{eq:IBP-vector},
\begin{equation}
    \langle \rho_t,F_u[f]\rangle_{L^2(\muS)}
    =
    \left\langle-\divS(\rho_tF_u),f\right\rangle_{L^2(\muS)}.
    \label{eq:FP-proof-drift}
\end{equation}
For the diffusion term,
\[
    \rho_t(\mathsf D:\Hess^{\mathrm S}f)
    =
    (\rho_t\mathsf D):\Hess^{\mathrm S}f.
\]
Applying the double tensor integration-by-parts identity \eqref{eq:tensor-double-IBP},
\begin{equation}
\begin{aligned}
    \int_{\mathsf{M}}
    (\rho_t\mathsf D):\Hess^{\mathrm S}f\,\dd\muS
    =
    \left\langle
    \divS\bigl(\divS(\rho_t\mathsf D)\bigr),
    f
    \right\rangle_{L^2(\muS)}.
\end{aligned}
\label{eq:FP-proof-diffusion}
\end{equation}
Substituting \eqref{eq:FP-proof-drift} and \eqref{eq:FP-proof-diffusion} into \eqref{eq:FP-proof-1},
\begin{equation*}
\begin{aligned}
    &\langle \rho_t,\mathcal A^u f\rangle_{L^2(\muS)}\\ 
    &\quad=
    \Big\langle
    -\divS(\rho_tF_u)
    +
    \frac12
    \divS\bigl(\divS(\rho_t\mathsf D)\bigr),
    f
    \Big\rangle_{L^2(\muS)}.
\end{aligned}
\end{equation*}
Hence
\[
    (\mathcal A^u)^*\rho_t
    =
    -\divS(\rho_tF_u)
    +
    \frac12
    \divS\bigl(\divS(\rho_t\mathsf D)\bigr).
\]
Combining this with the abstract forward equation \eqref{eq:abstract-FP} proves \eqref{eq:main-FP-double-div}.  Finally, inserting \eqref{eq:FP-D-repeat} and using the linearity of tensor divergence gives \eqref{eq:main-FP-HV}.  This is exactly the specialization of~\cite[Theorem~5]{LeeChirikjianFP2026} to the Sasaki state manifold and the horizontal--vertical diffusion tensor of the present controlled mechanical system.
\end{proof}

\subsection{Equivalence of the Stratonovich and It\^o forward equations}

The Stratonovich equation \eqref{eq:compact-Strat-for-FP} and the intrinsic It\^o equation \eqref{eq:intrinsic-Ito-state} describe the same diffusion process.  Indeed, the It\^o--Stratonovich conversion on $(\mathsf{M},\Sasaki)$ gives
\begin{equation}
    \widetilde X_u
    =
    X_u^{\mathrm{Strat}}
    +
    \frac12
    \sum_{i=1}^2
    \LCs_{B_i^q}B_i^q
    +
    \frac12
    \sum_{i=1}^2
    \LCs_{B_i^v}B_i^v.
    \label{eq:FP-Ito-Strat-conversion}
\end{equation}
Since the noise intensities are constant,
\begin{align*}
    \LCs_{B_i^q}B_i^q
    &=
    \sigma_q^2
    \LCs_{E_i^\Hor}E_i^\Hor
    =
    \sigma_q^2
    (\LC_{E_i}E_i)^\Hor,
    \\ 
    \LCs_{B_i^v}B_i^v
    &=
    \sigma_v^2
    \LCs_{E_i^\Ver}E_i^\Ver
    =
    0.
\end{align*}
Substitution into \eqref{eq:FP-Ito-Strat-conversion} yields
\[
    \widetilde X_u=F_u.
\]
Therefore the two SDE representations have the same generator and hence the same adjoint generator.  Consequently, the Stratonovich Fokker--Planck equation \eqref{eq:FP-Strat-compact} and the tensor-field It\^o equation \eqref{eq:main-FP-double-div} are two intrinsic representations of the same forward evolution.  This is the same equivalence mechanism emphasized in~\cite[Theorems~2--5]{LeeChirikjianFP2026}.

\subsection{Adjoints of the horizontal and vertical second-order operators}

Recall the backward operators
\[
    \Delta_{\Hor}f
    =
    \mathsf G_{\Hor}:\Hess^{\mathrm S}f,
    \qquad
    \Delta_{\Ver}f
    =
    \mathsf G_{\Ver}:\Hess^{\mathrm S}f.
\]
Their formal $L^2(\muS)$ adjoints are therefore
\begin{align}
    \Delta_{\Hor}^*\rho
    &:=
    \divS\!\left(
        \divS(\rho\mathsf G_{\Hor})
    \right),
    \label{eq:DeltaH-adjoint}\\
    \Delta_{\Ver}^*\rho
    &:=
    \divS\!\left(
        \divS(\rho\mathsf G_{\Ver})
    \right).
    \label{eq:DeltaV-adjoint}
\end{align}
Accordingly, \eqref{eq:main-FP-HV} may be written
\begin{equation*}
     {
    \frac{\partial \rho_t}{\partial t}
    =
    -\divS(\rho_tF_u)
    +
    \frac{\sigma_q^2}{2}\Delta_{\Hor}^*\rho_t
    +
    \frac{\sigma_v^2}{2}\Delta_{\Ver}^*\rho_t.
    }
\end{equation*}

\begin{remark}[Backward versus forward horizontal--vertical operators]
The operators $\Delta_{\Hor}$ and $\Delta_{\Ver}$ introduced earlier act on test functions and belong to the \emph{backward} generator.  The corresponding terms in the Fokker--Planck equation are their formal adjoints \eqref{eq:DeltaH-adjoint}--\eqref{eq:DeltaV-adjoint}.  We do not identify $\Delta_{\Hor}^*$ with $\Delta_{\Hor}$, or $\Delta_{\Ver}^*$ with $\Delta_{\Ver}$, without an additional self-adjointness argument.  The double-divergence tensor formula is the intrinsically correct forward operator and automatically includes derivatives of the horizontal and vertical tensor fields caused by the geometry.
\end{remark}

\subsection{Equal position and velocity noise: reduction to the Sasaki Laplacian}

The equal-intensity case admits an important simplification.  Suppose
\[
    \sigma_q=\sigma_v=\sigma.
\]
Then \eqref{eq:diffusion-tensor-HV} and \eqref{eq:cometric-splitting} give
\begin{equation*}
     {
    \mathsf D
    =
    \sigma^2(\Sasaki)^{-1}.
    }
\end{equation*}
Because $\LCs$ is metric compatible,
\[
    \LCs(\Sasaki)^{-1}=0.
\]
Hence
\begin{equation}
    \divS\!\left(
        \rho_t(\Sasaki)^{-1}
    \right)
    =
    \gradS \rho_t.
    \label{eq:div-p-cometric}
\end{equation}
Applying the divergence once more,
\begin{equation*}
    \divS\!\left[
        \divS\!\left(
            \rho_t(\Sasaki)^{-1}
        \right)
    \right]
    =
    \divS(\gradS \rho_t)
    =
    \Delta^{\mathrm S}\rho_t.
\end{equation*}
Therefore Theorem~\ref{thm:main-FP} reduces to
\begin{equation}
     {
    \frac{\partial \rho_t}{\partial t}
    =
    -\divS\!\left(
        \rho_t(v_q^\Hor+u^\Ver)
    \right)
    +
    \frac{\sigma^2}{2}
    \Delta^{\mathrm S}\rho_t.
    }
    \label{eq:FP-isotropic}
\end{equation}
This is the controlled advection--diffusion equation on the Sasaki manifold $(\mathsf{T}\Sph,\Sasaki)$: the first term transports probability along the deterministic controlled mechanical vector field, while the second term is isotropic Brownian diffusion on the full state manifold.

Equation~\eqref{eq:FP-isotropic} is also consistent with the backward generator \eqref{eq:generator-equal-noise}.  The Sasaki Laplace--Beltrami operator is formally self-adjoint with respect to $\muS$, so the same $\Delta^{\mathrm S}$ appears in both the backward generator and the forward equation.

\subsection{Probability flux and conservation of total probability}

The tensor-field form has a direct conservation-law interpretation.  Define the probability current
\begin{equation*}
     {
    J_t^u
    :=
    \rho_tF_u
    -
    \frac12
    \divS(\rho_t\mathsf D)
    \in\mathfrak X(\mathsf{M}).
    }
\end{equation*}
Then \eqref{eq:main-FP-double-div} is equivalent to
\begin{equation}
     {
    \frac{\partial \rho_t}{\partial t}
    +
    \divS J_t^u
    =
    0.
    }
    \label{eq:FP-continuity}
\end{equation}
Thus the Fokker--Planck equation is an intrinsic continuity equation on the tangent-bundle state manifold.  The first contribution
\[
    \rho_tF_u
\]
is the drift flux, whereas
\[
    -\frac12\divS(\rho_t\mathsf D)
\]
is the diffusion-induced flux.  This is exactly the geometric flux interpretation emphasized in~\cite[Remark~2]{LeeChirikjianFP2026}.

In the equal-intensity case,
\[
    \mathsf D=\sigma^2(\Sasaki)^{-1},
\]
so \eqref{eq:div-p-cometric} gives
\begin{equation*}
     {
    J_t^u
    =
    \rho_tF_u
    -
    \frac{\sigma^2}{2}\gradS \rho_t.
    }
\end{equation*}

To verify probability conservation, integrate \eqref{eq:FP-continuity} over $\mathsf{M}$:
\begin{equation*}
    \frac{\dd}{\dd t}
    \int_{\mathsf{M}} \rho_t\,\dd\muS
    =
    -\int_{\mathsf{M}}\divS J_t^u\,\dd\muS.
\end{equation*}
Since $\mathsf{M}=\mathsf{T}\Sph$ has no finite boundary, and assuming the current decays sufficiently rapidly in the noncompact fiber directions, the divergence theorem gives
\[
    \int_{\mathsf{M}}\divS J_t^u\,\dd\muS=0.
\]
Consequently,
\begin{equation*}
     {
    \frac{\dd}{\dd t}
    \int_{\mathsf{M}} \rho_t\,\dd\muS
    =
    0.
    }
\end{equation*}
Thus a normalized initial density remains normalized.

\subsection{Geometric interpretation of the anisotropic diffusion tensor}

The tensor
\[
    \mathsf D
    =
    \sigma_q^2\mathsf G_{\Hor}
    +
    \sigma_v^2\mathsf G_{\Ver}
\]
separates the covariance of the noise into the two geometrically distinguished subbundles
\[
    \mathsf{T}(\mathsf{T}\Sph)
    =
    \HH \mathsf{T}\Sph
    \oplus
    \VV \mathsf{T}\Sph.
\]
The horizontal tensor $\mathsf G_{\Hor}$ spreads probability in directions that move the configuration while parallel transporting the velocity state.  The vertical tensor $\mathsf G_{\Ver}$ spreads probability inside the tangent fibers and therefore represents uncertainty or random motion in the velocity variable.  The constants $\sigma_q^2$ and $\sigma_v^2$ are the corresponding diffusion strengths.

The inner divergence
\[
    \divS(\rho_t\mathsf D)
\]
is the intrinsic diffusion flux vector field generated by the density and the anisotropic covariance tensor.  The outer divergence measures the local accumulation or depletion caused by that flux.  In this form, curvature and the variation of the horizontal--vertical splitting enter through the Sasaki covariant derivative used in the two divergence operations rather than through coordinate Christoffel symbols.

\subsection{Frame independence and global character of the Fokker--Planck PDE}

The local Stratonovich representation uses a local orthonormal frame $\{E_1,E_2\}$ because $\Sph$ is not parallelizable.  Nevertheless, the tensor-field forward equation \eqref{eq:main-FP-double-div} is global.  The reason is that
\[
    F_u=v_q^\Hor+u^\Ver
\]
is a globally defined vector field whenever $u$ is a globally defined feedback section, and the diffusion tensor
\[
    \mathsf D
    =
    \sigma_q^2\mathsf G_{\Hor}
    +
    \sigma_v^2\mathsf G_{\Ver}
\]
is independent of the local orthonormal frame used to represent it.  Since $\divS$ is itself an intrinsic operator determined by $(\Sasaki,\LCs,\muS)$, every term in \eqref{eq:main-FP-double-div} is globally defined.

This is one of the main advantages of the tensor-field formulation in~\cite{LeeChirikjianFP2026}: although individual diffusion vector fields may only admit local frame representations, their covariance tensor and the resulting double-divergence forward operator are global geometric objects.

\subsection{Limiting cases of the Fokker--Planck equation}

The forward equation recovers the expected deterministic and partially stochastic limits.

\subsubsection{No stochastic noise}

If
\(
    \sigma_q=\sigma_v=0,
\)
then $\mathsf D=0$, and \eqref{eq:main-FP-double-div} becomes the intrinsic Liouville or continuity equation
\begin{equation*}
     {
    \frac{\partial \rho_t}{\partial t}
    =
    -\divS\!\left(
        \rho_t(v_q^\Hor+u^\Ver)
    \right).
    }
\end{equation*}

\subsubsection{Velocity noise only}

If
\(
    \sigma_q=0,
    \text{ and }
    \sigma_v>0,
\)
then
\[
    \mathsf D=\sigma_v^2\mathsf G_{\Ver},
\]
and
\begin{equation*}
     {
    \frac{\partial \rho_t}{\partial t}
    =
    -\divS(\rho_tF_u)
    +
    \frac{\sigma_v^2}{2}
    \divS\!\left(
        \divS(\rho_t\mathsf G_{\Ver})
    \right).
    }
\end{equation*}

\subsubsection{Position noise only}

If
\(
    \sigma_q>0,
    \text{ and }
    \sigma_v=0,
\)
then
\[
    \mathsf D=\sigma_q^2\mathsf G_{\Hor},
\]
and
\begin{equation*}
     {
    \frac{\partial \rho_t}{\partial t}
    =
    -\divS(\rho_tF_u)
    +
    \frac{\sigma_q^2}{2}
    \divS\!\left(
        \divS(\rho_t\mathsf G_{\Hor})
    \right).
    }
\end{equation*}

\subsubsection{Equal isotropic state noise}

If
\(
    \sigma_q=\sigma_v=\sigma,
\)
the forward equation is the Sasaki advection--diffusion equation \eqref{eq:FP-isotropic}.

\subsection{The intrinsic initial-value problem}

Given an initial probability density
\[
    \rho_0:\mathsf{M}\to[0,\infty)
\]
satisfying
\[
    \int_{\mathsf{M}}\rho_0\,\dd\muS=1,
\]
the forward problem associated with the controlled stochastic point mass is

\begin{equation}
\begin{aligned}
    &\frac{\partial \rho}{\partial t}
    ={}
    -\divS\!\left(
        \rho(v_q^\Hor+u^\Ver)
    \right)\\ 
    &\quad\quad\quad\quad\quad+\frac12
    \divS\!\Bigl[
        \divS\!\Bigl(
            \rho\bigl(\sigma_q^2\mathsf G_{\Hor} +\sigma_v^2\mathsf G_{\Ver}
            \bigr)
        \Bigr)
    \Bigr],\\[2mm]
    &\rho(0,\cdot)=\rho_0,\quad
    \rho(t,\cdot)\ge0,\quad \int_{\mathsf{M}}
    \rho(t,x)\,\dd\muS(x)=1.
\end{aligned}
\label{eq:FP-Cauchy-problem}
\end{equation}
Equation~\eqref{eq:FP-Cauchy-problem} is the intrinsic Fokker--Planck description associated with the intrinsic controlled stochastic dynamics developed in the preceding sections.

\subsection{Collected Final Form of the Fokker--Planck Equation}

Let $\rho_t$ denote the density of $X_t$ with respect to the Sasaki volume measure $\muS$, and define
\[
    \mathsf G_{\Hor}
    =
    \sum_{i=1}^2E_i^\Hor\otimes E_i^\Hor,
    \qquad
    \mathsf G_{\Ver}
    =
    \sum_{i=1}^2E_i^\Ver\otimes E_i^\Ver.
\]
The globally defined diffusion tensor is
\[
    \mathsf D
    =
    \sigma_q^2\mathsf G_{\Hor}
    +
    \sigma_v^2\mathsf G_{\Ver}.
\]
The associated Fokker--Planck PDE is
\begin{equation*}
     {
    \frac{\partial \rho_t}{\partial t}
    =
    -\divS\!\left(
        \rho_t(v_q^\Hor+u^\Ver)
    \right)
    +
    \frac12
    \divS\!\left(
        \divS(\rho_t\mathsf D)
    \right).
    }
\end{equation*}
For equal intensities $\sigma_q=\sigma_v=\sigma$,
\[
    \mathsf D=\sigma^2(\Sasaki)^{-1},
\]
and hence
\begin{equation*}
     {
    \frac{\partial \rho_t}{\partial t}
    =
    -\divS\!\left(
        \rho_t(v_q^\Hor+u^\Ver)
    \right)
    +
    \frac{\sigma^2}{2}
    \Delta^{\mathrm S}\rho_t.
    }
\end{equation*}

\section{Numerical Simulations}
\label{sec:Numerical_simulations}
\begin{figure*}[!t]
    \centering
    \includegraphics[scale = 0.5]{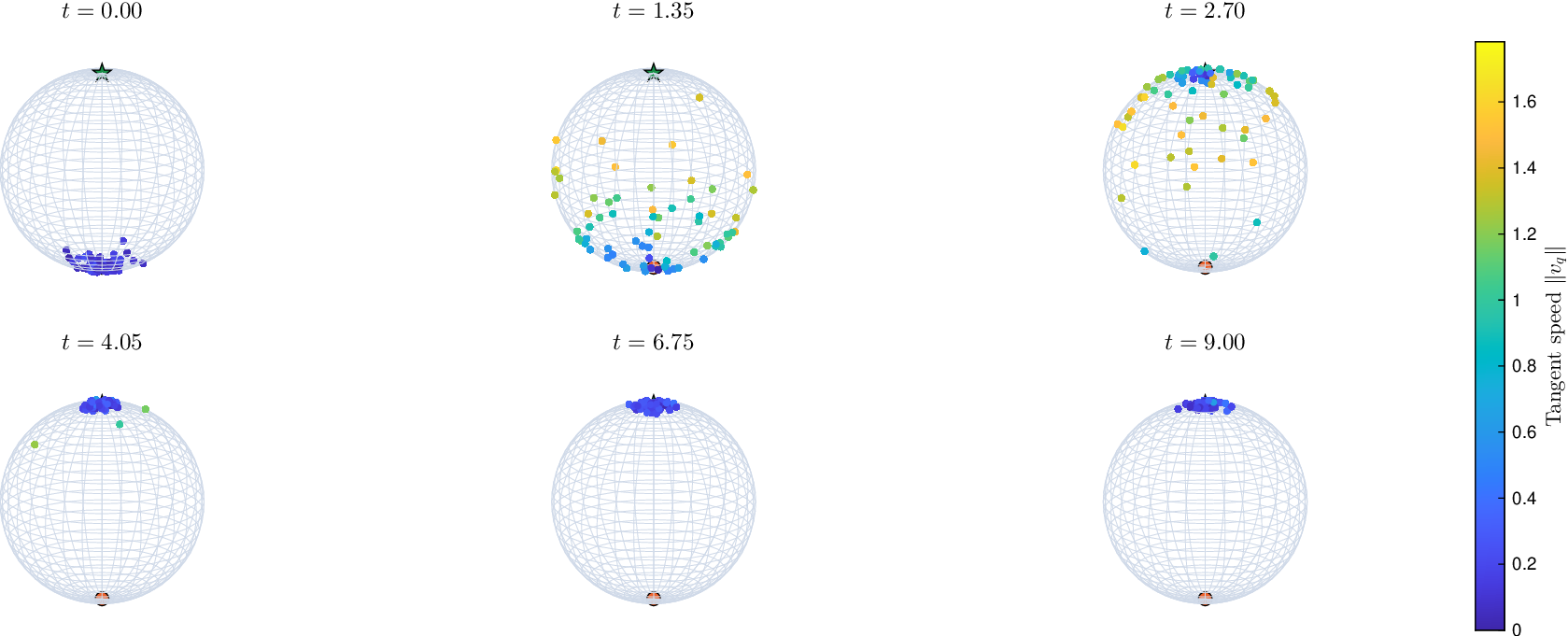}
\caption{Time snapshots of the stochastic particle
ensemble on $\mathbb S^2$ at six different times. Each point represents the position of one particle,
with its color indicating the magnitude of the corresponding
tangent velocity. The controlled ensemble evolves from an
initial distribution concentrated near $q_{\mathrm I}$
toward a neighborhood of $q_\star$, while stochastic
perturbations maintain a nonzero distributional spread.}
\label{fig:snapshots__particles}
\end{figure*}

We illustrate the controlled stochastic dynamics developed in this paper by
simulating an ensemble of $N$ independent point masses on $\Sph$ with
independent horizontal and vertical disturbances. The objective is to drive
an initial ensemble concentrated near $0_{q_{\mathrm I}}\in \mathsf{T}\Sph$ toward
a neighborhood of $0_{q_\star}\in \mathsf{T}\Sph$. The state distribution thereby
evolves according to the Fokker--Planck equation derived above; the
experiment approximates this evolution through sample trajectories rather
than a direct numerical solution of the four-dimensional PDE.

\emph{Feedback controller and initial distribution:}
We employ the intrinsic potential-based feedback
\begin{equation*}
\begin{aligned}
    U(q)&=k_p\bigl(1-\cos d_{\Sph}(q,q_\star)\bigr),\\
    u(v_q)&=-\operatorname{grad}_{\Grnd}U(q)-k_d v_q,
\end{aligned}
\end{equation*}
where $d_{\Sph}$ is the geodesic distance on the unit sphere and
$k_p,k_d>0$ are the attraction and damping gains, respectively.
The potential is globally smooth, although its gradient vanishes at the
antipode of $q_\star$ as well as at the target. In the
numerical realization, we take antipodal points $q_{\mathrm I}=(0,0,-1)$ and
$q_\star=(0,0,1)$ as the initial and target configurations, respectively. Initial conditions are sampled
independently according to
\begin{equation*}
\begin{aligned}
    q_0^{(j)}&=\operatorname{Exp}_{q_{\mathrm I}}
       (\varepsilon_q Z_0^{(j)}),\\
    v_{q_0^{(j)}}^{(j)}&=\varepsilon_v\widetilde Z_0^{(j)},
    \qquad j=1,\ldots,N,
\end{aligned}
\end{equation*}
where $Z_0^{(j)}$ is an isotropic standard Gaussian tangent vector at
$q_{\mathrm I}$ and $\widetilde Z_0^{(j)}$, conditioned on
$q_0^{(j)}$, is an independent isotropic standard Gaussian vector in
$\mathsf{T}_{q_0^{(j)}}\Sph$. The map $\operatorname{Exp}_{q_{\mathrm I}}$ denotes the Riemannian
exponential map at $q_{\mathrm I}$, which maps a tangent vector
in $\mathsf{T}_{q_{\mathrm I}}\Sph$ to a point on $\Sph$
along the corresponding geodesic. Thus the initial joint distribution is
concentrated about $0_{q_{\mathrm I}}$, with position and velocity
spreads $\varepsilon_q$ and $\varepsilon_v$, respectively.

\emph{Geometric stochastic integration:}
Let $t_k=k\Delta t$ and write
$X_k^{(j)}=v_{q_k^{(j)}}^{(j)}$. For each particle and time step,
independent isotropic Gaussian tangent vectors
$Z_{q,k}^{(j)},Z_{v,k}^{(j)}\in \mathsf{T}_{q_k^{(j)}}\Sph$ are sampled. The
intrinsic It\^o state equation is approximated by the geometric
Euler-type update
\begin{equation*}
\begin{aligned}
    \delta q_k^{(j)}
      &=v_{q_k^{(j)}}^{(j)}\Delta t
        +\sigma_q\sqrt{\Delta t}\,Z_{q,k}^{(j)},\\
    \widetilde v_{q_k^{(j)},k+1}^{(j)}
      &=v_{q_k^{(j)}}^{(j)}+u(X_k^{(j)})\Delta t
        +\sigma_v\sqrt{\Delta t}\,Z_{v,k}^{(j)},\\
    q_{k+1}^{(j)}
      &=\operatorname{Exp}_{q_k^{(j)}}(\delta q_k^{(j)}),\\
    v_{q_{k+1}^{(j)}}^{(j)}
      &=\mathcal P_{q_k^{(j)}\to q_{k+1}^{(j)}}
        (\widetilde v_{q_k^{(j)},k+1}^{(j)}),
\end{aligned}
\end{equation*}
where $\operatorname{Exp}_{q_k^{(j)}}$ is the Riemannian exponential map at ${q_k^{(j)}}$ and $\mathcal P_{q_k^{(j)}\to q_{k+1}^{(j)}}$ denotes
the parallel transport map along the shortest geodesic
connecting $q_k^{(j)}$ to $q_{k+1}^{(j)}$ on $\Sph$. For numerical evaluation and visualization, these intrinsic
operations are implemented using the standard three-dimensional unit-vector
representation of $\Sph$; Gaussian vectors are orthogonally projected
onto the corresponding tangent planes. After each update, numerical
roundoff is removed by normalizing the configuration and projecting the
velocity onto its tangent space. Thus the algorithm maintains the
constraints $q_k^{(j)}\in\Sph$ and
$v_{q_k^{(j)}}^{(j)}\in \mathsf{T}_{q_k^{(j)}}\Sph$.

\emph{Parameters and diagnostics:}
The experiment simulates
$N=70$ particles with a time horizon of $T=9$, a time step
of $\Delta t=0.0075$, control gains $k_p=7$ and $k_d=4.5$,
noise intensities $\sigma_q=0.12$ and $\sigma_v=0.32$,
and initial position and velocity spreads
$\varepsilon_q=0.15$ and $\varepsilon_v=0.10$, respectively.
Figure~\ref{fig:snapshots__particles} shows the time snapshots of the stochastic particle
ensemble on $\mathbb S^2$ at six different times. The controlled ensemble evolves from an
initial distribution concentrated near $q_{\mathrm I}$
toward a neighborhood of $q_\star$, while stochastic
perturbations maintain a nonzero distributional spread. Figure~\ref{fig:sample_path_trajectories} shows the sample path trajectories of each particle in the ensemble.

\section{Concluding Remarks}
\label{sec:Conclusion}
This paper developed an intrinsic geometric stochastic model
of a fully actuated point mass on $\Sph$, with state manifold
$\mathsf{T}\Sph$. Using the Levi--Civita horizontal--vertical
splitting and Sasaki geometry, we formulated controlled
Stratonovich and It\^o dynamics with independent isotropic
position and velocity noise. The associated infinitesimal
generator was derived in horizontal--vertical Sasaki-Hessian
form, reducing to the Sasaki Laplace--Beltrami diffusion
operator for equal noise intensities. The corresponding Fokker--Planck equation was derived
intrinsically with respect to the Sasaki volume measure.
A globally defined horizontal--vertical diffusion tensor
yielded a coordinate-free double-divergence formulation,
which specializes to the Sasaki advection--diffusion equation
for equal noise intensities. Numerical simulations illustrated
the evolution of a controlled stochastic particle ensemble
toward an antipodal target configuration on $\Sph$.

\textit{Generative AI Disclosure:} ChatGPT~\cite{openai_chatgpt_2026} has been used to improve the syntax and grammar of some paragraphs.

\bibliographystyle{IEEEtran}
\bibliography{references}

\end{document}